\documentclass[reqno]{amsart}
\usepackage{amsmath,amssymb,stmaryrd}
\begin{document}
\newtheorem{thm}{Theorem}[section]
\newtheorem{lem}{Lemma}[section]
\newtheorem{prop}{Proposition}[section]
\newtheorem{coro}{Corollary}[section]
\newtheorem{defi}{Definition}[section]
\def \eps {\epsilon}
\def \rr {\mathbb{R}}
\def \rn {\mathbb{R}^n}
\def \rnm {\mathbb{R}^n_-}
\def \we {w_\eps}
\def \ze {z_\eps}
\def \he {h_\eps}
\def \xune {x_{1,\eps}}
\def \pe {p_\eps}
\def \tge {\tilde{g}_\eps}
\def \bge {\bar{g}_\eps}
\def \lae {\lambda_\eps}
\def \xe {x_\eps}
\def \ye {y_\eps}
\def \fe {f_\eps}
\def \ae {a_\eps}
\def \ue {u_\eps}
\def \me {\mu_\eps}
\def \ke {k_\eps}
\def \be {\beta_\eps}
\def \ve {v_\eps}
\def \tv {\tilde{v}}
\def \tve {\tilde{v}_\eps}
\def \twe {\tilde{w}_\eps}
\def \tw {\tilde{w}}
\def \bwe {\overline{w}_\eps}
\def \bw {\overline{w}}
\def \bve {\overline{v}_\eps}
\def \tye {\tilde{y}_\eps}
\def \tue {\tilde{u}_\eps}
\def \txe {\tilde{x}_\eps}
\def \tze {\tilde{z}_\eps}
\def \bze {\bar{z}_\eps}
\def \crit {2^{\star}}
\def \huno {H_{1,0}^2(\Omega)}
\def \hunrn {H_{1,0}^2(\rn)}
\def \hunrnm {H_{1,0}^2(\rnm)}
\def \ds {\displaystyle}
\def \beq {\begin{eqnarray*}}
\def \eeq {\end{eqnarray*}}
\def \beqn {\begin{eqnarray}}
\def \eeqn {\end{eqnarray}}
\def \bequa {\begin{equation}}
\def \eequa {\end{equation}}

\title[Best Hardy-Sobolev constant]{The effect of curvature on the best
constant in the Hardy-Sobolev inequalities}
\author{N. Ghoussoub}
\address{Nassif Ghoussoub, Department of Mathematics, University of
British Columbia, Vancouver, Canada}
\thanks{Research partially supported by the Natural Sciences and
Engineering Research Council of Canada. The first named author gratefully
acknowledges the hospitality and support of the Universit\'e de Nice where
this work was initiated.}
\email{nassif@math.ubc.ca}
\author{F. Robert}
\address{Fr\'ed\'eric Robert, Laboratoire J.A.Dieudonn\'e, Universit\'e de
Nice Sophia-Antipolis,
Parc Valrose, 06108 Nice cedex 2, France}
\email{frobert@math.unice.fr}
\thanks{The second named author gratefully acknowledges the hospitality
and support of the University of British Columbia where this work was
completed.}
\date{February 25th 2005}

\begin{abstract} We address the question of attainability of  the best
constant in the following Hardy-Sobolev inequality on a smooth domain
$\Omega$ of $\rn$:
$$ \mu_s (\Omega) := \inf \left \{ \int_{\Omega}| \nabla u|^2
dx;\, u \in \huno \hbox{ and }  \int_{\Omega} \frac {|u|^{\crit}}{|x|^s}\,
dx =1\right\}$$
when $0<s<2$, $\crit :=2^*(s)=\frac{2(n-s)}{n-2}$, and when $0$ is on the 
boundary
$\partial \Omega$. This question is closely related to the geometry of
$\partial\Omega$, as we extend here  the main result
obtained in \cite{gk} by proving that at least in dimension $n\geq 4$, the 
negativity of the
mean curvature   of $\partial \Omega $ at $0$ is sufficient to ensure the
attainability of $\mu_{s}(\Omega)$. Key ingredients in our proof are the 
identification of symmetries enjoyed by the extremal functions 
correrresponding to the best constant in half-space, as well as a fine 
analysis of the asymptotic behaviour of appropriate minimizing sequences. 
The result holds true also  in dimension 3 but  the more involved  proof 
will
be dealt with in a forthcoming paper \cite{gr2}.
\end{abstract}

\maketitle

\section{Introduction} Let $\Omega$ be a smooth domain of $\rn$, $n\geq 3$
and  denote by $\huno$ the completion of $C_c^\infty(\Omega)$, the set of
smooth functions compactly supported in $\Omega$, for the norm
$\Vert u\Vert_{\huno}=\sqrt{\int_{\Omega}|\nabla u|^2\, dx}$. The
Hardy-Sobolev inequality (\cite{ca}, \cite{cw}, \cite{gy}) asserts that
for $s\in [0,2]$ and for $\crit:=2^*(s)=\frac{2(n-s)}{n-2}$, there exists 
$C>0$
such that for all $u\in \hunrn$,
\bequa\label{ineq:HS:rn}
\left(\int_{\rn}\frac{|u|^{\crit}}{|x|^s}\,
dx\right)^{\frac{2}{\crit}}\leq C \int_{\rn}|\nabla u|^2\, dx.
\eequa
We define
\bequa\label{def:mus}
\mu_{s}(\Omega)=\inf\left\{\frac{\int_{\Omega}|\nabla u|^2\,
dx}{\left(\int_{\Omega}\frac{|u|^{\crit}}{|x|^s}\,
dx\right)^{\frac{2}{\crit}}}; u\in \huno\setminus\{0\}\right\},
\eequa
and we consider the corresponding ground state solutions in $\huno\cap
C^1(\overline{\Omega})$ for
\bequa\label{def:equa}
\left\{\begin{array}{ll}
\Delta u=\frac{u^{\crit-1}}{|x|^s}& \hbox{ in }{\mathcal D}'(\Omega)\\
u>0&\hbox{ in }\Omega\\
u=0&\hbox{ on }\partial\Omega.
\end{array}\right.
\eequa
where here and throughout the paper, $\Delta=-\sum_i\partial_{ii}$ is the
Laplacian with minus sign convention.

It is well known (see for instance \cite{st}) that in the non-singular
case i.e., when $s=0$, we have $\mu_s(\Omega) =\mu_s(\rn)$ for any domain
$\Omega$ and that $\mu_s(\Omega)$ is never attained unless ${\rm
cap}(\rn\setminus \Omega) =0$. In this situation, the limiting spaces
after blow-up of solutions of (\ref{def:equa}) is $\rn$. It was shown in
\cite{gy} that the same result holds true for any $0<s<2$ as long as $0$
belongs to the interior of a domain.

However, the fact that things may be different when $0\in \partial \Omega$
first emerged in a paper by Egnell \cite{eg} where  he considers open
cones of the form ${C}=\{x\in \rn; x=r\theta, \theta\in D \ {\rm and}\
r>0\}$ where the base $D$ is a connected domain of the unit sphere
$S^{n-1}$ of $\rn$.  Egnell showed that $\mu_{s}(C)$ is then attained for
$0<s<2$ even when $\bar {C} \neq \rn$. This obviously applies to
  a half-space
$\rnm=\{x\in\rnm/\, x_1<0\}$, where $x_1$ denotes the first coordinate of
a generic point $x\in\rn$ in the canonical basis of $\rn$.

Half-spaces containing $0$ on their boundary were identified in \cite{gk} 
as the limiting spaces after blow-up
in the case where $\partial \Omega$ is smooth at $0$, and the curvature of 
the
boundary at $0$ then gets to play an important role. In our context, we 
specify
the orientation of $\partial\Omega$ in such a way that the normal vectors
of $\partial \Omega$ are pointing outward from the domain $\Omega$.
It was shown in \cite{gk} that in dimension $n\geq 4$, the negativity of
all principal curvatures at $0$ --which is essentially a condition of {\it
``strict concavity''}  at 0-- leads to attainability of the best constant
for problems with Dirichlet boundary conditions, while the Neumann
problems required the positivity of the mean curvature at $0$. On the
other hand, standard Pohozaev type arguments show non-attainability in the
cases where $\Omega$ is convex or star-shaped at $0$.
 
In this paper, we improve and complete the results in \cite{gk} in a 
substantial  way  by showing that   for the  best constant to be achieved, 
it is sufficient that the mean curvature be negative. This is now quite 
similar but dual  to the case with Neumann boundary conditions which 
requires the  mean curvature to be positive.

More precisely, assume that the principal curvatures
$\alpha_1,...,\alpha_{n-1}$ of $\partial \Omega$ at $0$ are finite. The
oriented boundary $\partial\Omega$ near the origin can then be represented
(up to rotating the coordinates if necessary) by $x_1 = \varphi_0 (x') =
-\frac{1}{2} \sum_{i=2}^{n}\alpha_{i-1} x^2_i +o(|x'|^2)$, where $x'=(x_2,
...,x_{n}) \in B_\delta(0) \cap \{ x_1 =0 \}$ for some $\delta >0$ and 
where
$B_\delta(0)$ is the ball in $\rn$ centered at $0$ with radius $\delta$.
If one assumes  the principal curvatures at $0$ to be negative, that is if
$$\max_{1 \leq i \leq n-1}\alpha_i <0,$$
then the sectional curvature at $0$ is negative and therefore $\partial
\Omega$ --viewed as an $(n-1)$-Riemannian submanifold of $\rn$-- is
strictly convex at $0$ (see for instance \cite{ghl}). The latter property
means that there exists a neighborhood $U$ of $0$ in $\partial \Omega$,
such that the whole of $U$ lies on one side of a hyperplane $H$ that is
tangent to $\partial \Omega$ at $0$ and $U\cap H=\{0\}$, and so does the
complementary $\rn\setminus\Omega$, at least locally. The above curvature
condition then amounts to a notion of strict local convexity of $\rn
\setminus \Omega$ at $0$. Our main result below shows that at least for
dimension greater than $4$, it is sufficient to assume that
$$\sum_{1 \leq i \leq n-1}\alpha_i < 0.$$
\begin{thm}\label{th:intro}
Let $\Omega$ be a smooth bounded oriented domain of $\rn$ where $n\geq 4$, 
such
that $0\in\partial\Omega$ and assume $s\in (0,2)$.  If the mean curvature 
of $\partial\Omega$ at $0$ is negative, then the infimum $\mu_{s}(\Omega)$ 
in (\ref{def:mus}) is achieved. In addition, the set of minimizers of 
(\ref{def:mus}) is pre-compact in the $\huno-$topology.
\end{thm}
The first difficulty we have to face here is that   the extremals  for 
(\ref{def:mus}) when $\Omega=\rnm$ are not known explicitely, and our 
first result below --proved in section 2-- is the identification of 
certain symmetries enjoyed by these extremals --and actually all positive 
solutions-- on half-space.

  \begin{thm}\label{prop:sym}
Let $n\geq 3$, $s\in (0,2)$ and consider $u\in C^2(\rnm)\cap 
C^1(\overline{\rnm})$ such that
\bequa\label{sys:sym}\left\{\begin{array}{ll}
\Delta u=\frac{u^{\crit-1}}{|x|^s}& \hbox{ in }\rnm\\
\quad u>0 & \hbox{ in }\rnm\\
\quad u=0 & \hbox{ on }\partial\rnm,
\end{array}
\right.\eequa
where $\crit=\frac{2(n-s)}{n-2}$. Assume that for some $C>0$, $u(x)\leq 
C(1+|x|)^{1-n}$ for all $x\in\rnm$. Then we have that $u\circ\sigma=u$ for 
all isometry of $\rn$ such that $\sigma(\rnm)=\rnm$. In particular, there 
exists $v\in C^2(\rr_-^\star\times \rr)\cap C^1(\rr_-\times \rr)$ such 
that for all $x_1<0$ and all $x'\in\rr^{n-1}$, we have that 
$u(x_1,x')=v(x_1,|x'|)$.
\end{thm}
 
The  attainability  result is then obtained by combining this new 
information with a fine study of the asymptotic behaviour of
solutions to the corresponding subcritical pde's. They can eventually 
develop a singularity at zero as we approach the critical exponent 
$2^*(s)$, and for that we proceed to completely describe the way they may 
blow up, which makes for an interesting analysis in its own right.\\
  Indeed, assume $\Omega$ is a smooth bounded  domain of $\rn$ such that 
$0\in\partial\Omega$ and consider for any $\eps\in (0,\crit-2)$,
  the infimum
$$\mu_{s}^\eps(\Omega):=\inf_{u\in\huno\setminus\{0\}}\frac{\int_\Omega
|\nabla u|^2\,
dx}{\left(\int_\Omega\frac{|u|^{\crit-\eps}}{|x|^s}\,
dx\right)^{\frac{2}{\crit-\eps}}},$$
which is achieved by a function $\ue\in\huno$, $\ue> 0$ in $\Omega$ in
$C^1(\overline{\Omega})\cap C^2(\overline{\Omega}\setminus\{0\})$
that satisfies the system
$$\left\{\begin{array}{ll|}
  \Delta\ue=\frac{\ue^{\crit-1-\eps}}{|x|^s}  &\hbox{ in }{\mathcal
D}'(\Omega)\\
\quad   \ue>0 &  \hbox{ in }\Omega\\
\int_\Omega\frac{|u_\eps|^{\crit-\eps}}{|x|^s}\,
dx=(\mu_{s}^\eps(\Omega))^{\frac{\crit-\eps}{\crit-2-\eps}}.&
\end{array}\right.$$
  The bulk of the paper (beyond section 2)  consists of  proving the 
following estimate.
\begin{thm}\label{prop:poho.0} Let $\Omega$ be a smooth bounded oriented 
domain of $\rn$ where $n\geq 4$,
and assuming that $u_\epsilon$ converges weakly to zero (i.e. when blow-up 
occurs),
  then there exists $v$ solution for (\ref{sys:sym}) such that
  \[
  \int_{\rnm}|\nabla v|^2\, 
dx=\mu_{s}(\Omega)^{\frac{\crit}{\crit-2}}=\mu_{s}(\rnm)^{\frac{\crit}{\crit-2}},
\]
while  -modulo passing to a subsequence- we have
$$\lim_{\eps\to 0}\eps \, (\max_\Omega 
u_\eps)^{\frac{2}{n-2}}=\frac{(n-s)\int_{\partial\rnm}|x|^2|\nabla v|^2\, 
dx}{n(n-2)^2\mu_{s}(\rnm)^{\frac{n-s}{2-s}}}\cdot H(0).$$
where  $H(0)$ is the mean curvature of the oriented boundary 
$\partial\Omega$ at $0$.
\end{thm}
  These techniques actually allow us to
prove the following existence theorem. We shall say that a function is
in $C^1(\overline{\Omega})$ if it can be extended to a $C^1-$function in a
neighborhood of $\Omega$.

\begin{thm}\label{th:eq}
Let $\Omega$ be a smooth bounded oriented domain of $\rn$ where $n\geq 4$, 
such that $0\in\partial\Omega$. Assume $s\in (0,2)$ and consider $a\in 
C^1(\overline{\Omega})$ such that the operator $\Delta+a$ is coercive in 
$\Omega$. If the mean curvature of $\partial\Omega$ at $0$ is negative, 
then there exists a solution  $u\in \huno\cap C^1(\overline{\Omega})$ for
$$\left\{\begin{array}{ll}
\Delta u+au=\frac{u^{\crit-1}}{|x|^s}& \hbox{ in }{\mathcal D}'(\Omega)\\
u>0&\hbox{ in }\Omega\\
u=0&\hbox{ on }\partial\Omega.
\end{array}\right.$$
\end{thm}

  The study of blow-up solutions in certain nonlinear elliptic equations
was initiated by Atkinson-Peletier \cite{ap} (see also Br\'ezis-Peletier
\cite{bp}). In the Riemannian context, such asymptotics were first 
studied by Schoen \cite{s} and Hebey-Vaugon \cite{hv1}. The techniques of
blow-up have been developed in a general context by Druet, Hebey and the
second author \cite{dhr}. They turned out to be very powerful
tools for the study of best constant problems in Sobolev inequalities, see
for instance Druet \cite{d1}, Hebey-Vaugon \cite{hv1}, \cite{hv2} and
Robert \cite{r2}). We also mention the work of  Han \cite{ha}, Hebey \cite{he}, Druet-Robert
\cite{dr} and Robert \cite{r1}) on the asymptotics for
solutions to nonlinear pde's, the $3-$dimensional conjecture of
Br\'ezis solved by Druet \cite{d2} and the intricate compactness issues in
the Riemannian context (see for instance Schoen \cite{s} and Druet
\cite{d3}).
 
In a forthcoming paper \cite{gr1}, we shall establish a more refined
compactness result which yields an infinite number of sign
changing solutions for (\ref{def:equa}). In another forthcoming article
\cite{gr2}, we tackle similar questions for various critical equations
involving a whole affine subspace of singularities on the boundary.

\section{Symmetry of the positive solutions to the limit 
equation}\label{sec:sym}

This section is devoted to the proof of Theorem \ref{prop:sym}, that is 
the symmetry property for the positive solutions to the limit equation on 
$\rnm$.
 For that, we consider $u\in C^2(\rnm)\cap C^1(\overline{\rnm})$ that 
verifies the system (\ref{sys:sym}) while verifying for some $C>0$ the  
bound
\bequa\label{ineq:sym}
u(x)\leq \frac{C}{(1+|x|)^{n-1}}
\eequa
for all $x\in\rnm$. Denoting by $\vec{e}_1$ the first vector of the canonical 
basis of $\rn$,  we consider the open ball $D:=B_{1/2}\left(-\frac{1}{2}\vec{e}_1\right)$
and define
\bequa
v(x):=|x|^{2-n}u\left(\vec{e}_1+\frac{x}{|x|^2}\right)
\eequa
for all $x\in \overline{D}\setminus\{0\}$ and $v(0)=0$. 
Clearly, this is well-defined.

\medskip\noindent{\bf Step 2.1:} We claim that
\bequa\label{ineq:hopf}
v\in C^2(D)\cap C^1(\overline{D})\hbox{ and }\frac{\partial v}{\partial 
\nu}<0\hbox{ on }\partial D
\eequa
where $\partial/\partial\nu$ denotes the outward normal derivative.
\begin{proof} It follows from the assumptions on $u$ that $v\in C^2(D)\cap 
C^1(\overline{D}\setminus\{0\})$. Moreover, $v(x)>0$ for all $x\in D$ and 
$v(x)=0$ for all $x\in\partial D\setminus\{0\}$. It follows from  
 (\ref{ineq:sym}) that there exists $C>0$ such that
\bequa\label{ineq:v:1}
v(x)\leq C|x|
\eequa
for all $x\in \overline{D}\setminus \{0\}$. Since $v(0)=0$, we have that 
$v\in C^0(\overline{D})$. The function $v$ verifies the equation
\bequa\label{eq:v}
\Delta 
v=\frac{v^{\crit-1}}{|x+|x|^2\vec{e}_1|^s}=\frac{v^{\crit-1}}{|x|^s\left|x+\vec{e}_1\right|^s}
\eequa
in $D$. Since $-\vec{e}_1\in \partial D\setminus\{0\}$ and $v\in 
C^1(\overline{D}\setminus\{0\})\cap C^0(\overline{D})$, there exists $C>0$ 
such that
\bequa\label{ineq:v:2}
v(x)\leq C |x+\vec{e}_1|
\eequa
for all $x\in \overline{D}$. It then follows from (\ref{ineq:v:1}), 
(\ref{eq:v}), (\ref{ineq:v:2}) and standard elliptic theory that $v\in 
C^1(\overline{D})$. Since $v>0$ in $D$, it follows from Hopf's Lemma that 
$\frac{\partial v}{\partial \nu}<0$ on $\partial D$.\end{proof}

We prove the symmetry of $u$ by proving a symmetry property of $v$, which 
is defined on a ball. Our proof uses the moving plane method. We take 
largely inspiration in \cite{gnn} and \cite{cgs}. Classically, for any 
$\mu\geq 0$ and any $x=(x',x_n)\in \rn$ ($x'\in \rr^{n-1}$ and 
$x_n\in\rr$), we let
$$x_\mu=(x',2\mu-x_n)\hbox{ and }D_\mu=\{x\in D/\, x_\mu\in D\}.$$
It follows from Hopf's Lemma (See (\ref{ineq:hopf}))
that there exists $\epsilon_0>0$ such that for any $\mu\in 
(\frac{1}{2}-\epsilon_0,\frac{1}{2})$, we have that $D_\mu\neq\emptyset$ 
and $v(x)\geq v(x_\mu)$ for all $x\in D_\mu$ such that $x_n\leq\mu$. We 
let $\mu\geq 0$. We say that $(P_\mu)$ holds if:
\[
\hbox{ $D_{\mu}\neq \emptyset$ and
$v(x)\geq v(x_{\mu})$
for all $x\in D_{\mu}$ such that $x_n\leq\mu$.}
\]
We let
\bequa\label{def:lambda}
\lambda:=\min\left\{\mu\geq 0;\,  (P_{\nu})\hbox{ holds for all }\nu\in 
\left(\mu,\frac{1}{2}\right) \right\}.
\eequa

\medskip\noindent{\bf Step 2.2:} We claim that $\lambda=0$.
\begin{proof} We proceed by contradiction and assume that $\lambda>0$. We 
then get that $D_{\lambda}\neq\emptyset$ and that $(P_{\lambda})$ holds. 
We let
$$w(x):=v(x)-v(x_{\lambda})$$
for all $x\in D_{\lambda}\cap\{x_n<\lambda\}$. Since $(P_{\lambda})$ 
holds, we have that $w(x)\geq 0$ for all $x\in 
D_{\lambda}\cap\{x_n<\lambda\}$. With the equation (\ref{eq:v}) of $v$ and 
$(P_{\lambda})$, we get that
\beq
\Delta w&=& \frac{v(x)^{\crit-1}}{|x+|x|^2\vec{e}_1|^s}- 
\frac{v(x_{\lambda})^{\crit-1}}{|x_{\lambda}+|x_{\lambda}|^2\vec{e}_1|^s}\\
&\geq &v(x_{\lambda})^{\crit-1}\left(\frac{1}{|x+|x|^2\vec{e}_1|^s}- 
\frac{1}{|x_{\lambda}+|x_{\lambda}|^2\vec{e}_1|^s}\right)
\eeq
for all $x\in D_{\lambda}\cap\{x_n<\lambda\}$. With straightforward 
computations, we have that
\beq
&& |x_{\lambda}|^2-|x|^2=4\lambda(\lambda-x_n)\\
&& 
|x_{\lambda}+|x_{\lambda}|^2\vec{e}_1|^2-|x+|x|^2\vec{e}_1|^2=(|x_{\lambda}|^2-|x|^2)\left(1+|x_{\lambda}|^2+|x|^2+2x_1)\right)
\eeq
for all $x\in\rn$. It follows that $\Delta w(x)>0$ for all $x\in 
D_{\lambda}\cap\{x_n<\lambda\}$. Note that we have used that $\lambda>0$. 
It then follows from Hopf's Lemma and the strong comparison principle that
\bequa\label{ppty:w}
w>0\hbox{ in }D_\lambda\cap\{x_n<\lambda\}\hbox{ and }\frac{\partial 
w}{\partial \nu}<0\hbox{ on }D_\lambda\cap\{x_n=\lambda\}.
\eequa
By definition, there exists a sequence 
$(\lambda_i)_{i\in\mathbb{N}}\in\rr$ and a sequence 
$(x^i)_{i\in\mathbb{N}}\in D$ such that $\lambda_i<\lambda$, $x^i\in 
D_{\lambda_i}$, $(x^i)_n<\lambda_i$, $\lim_{i\to 
+\infty}\lambda_i=\lambda$ and
\bequa\label{ineq:sym:2}
v(x^i)<v((x^i)_{\lambda_i})
\eequa
for all $i\in\mathbb{N}$. Up to extraction a subsequence, we assume that 
there exists $x\in \overline{D}_{\lambda}\cap\{x_n\leq \lambda\}$ such 
that $\lim_{i\to +\infty}x^i=x$ with $x_n\leq \lambda$. Passing to the 
limit $i\to +\infty$ in (\ref{ineq:sym:2}), we get that $v(x)\leq 
v(x_{\lambda})$. It follows from this last inequality and (\ref{ppty:w}) 
that $v(x)-v(x_\lambda)=w(x)=0$, and then $x\in 
\partial(D_{\lambda}\cap\{x_n<\lambda\})$.

\smallskip\noindent{\it Case 1:} If $x\in\partial D$. Then 
$v(x_{\lambda})=0$ and $x_{\lambda}\in\partial D$. Since $D$ is a ball and 
$\lambda>0$, we get that $x=x_{\lambda}\in\partial D$. Since $v$ is $C^1$, 
we get that there exists $\tau_i\in ((x^i)_n,2\lambda_i-(x^i)_n)$ such 
that
$$v(x^i)-v((x^i)_{\lambda_i})=\partial_n v((x')^i,\tau_i)\times 
2((x^i)_n-\lambda_i)$$
Letting $i\to +\infty$, using that $(x^i)_n<\lambda_i$ and 
(\ref{ineq:sym:2}), we get that $\partial_n v(x)\geq 0$. On the other 
hand, we have that
$$\partial_n v(x)=\frac{\partial v}{\partial\nu}(x)\cdot 
(\nu(x)|\vec{e}_n)=\frac{\lambda}{|x+\vec{e}_1/2|}\frac{\partial 
v}{\partial\nu}(x)<0.$$
A contradiction with (\ref{ineq:hopf}).

\smallskip\noindent{\it Case 2:} If $x\in D$. Since $v(x_{\lambda})=v(x)$, 
we then get that $x_{\lambda}\in D$. Since $x\in 
\partial(D_{\lambda}\cap\{x_n<\lambda\})$, we then get that $x\in 
D\cap\{x_n=\lambda\}$. With the same argument as in the preceding step, we 
get that $\partial_n v(x)\geq 0$. On the other hand, with (\ref{ppty:w}), we get that $2\partial_n v(x)=\partial_n w(x)<0$. A contradiction.

\smallskip\noindent In all the cases, we have obtained a contradiction. 
This proves that $\lambda=0$. \end{proof}

\medskip\noindent{\bf Step 2.3:} Here goes the final argument. Since 
$\lambda=0$, it follows from the definition (\ref{def:lambda}) of 
$\lambda$ that $v(x',x_n)\geq v(x',-x_n)$ for all $x\in D$ such that 
$x_n\leq 0$. With the same technique, we get the reverse inequality, and 
then, we get that
$$v(x',x_n)=v(x',-x_n)$$
for all $x=(x',x_n)\in D$. In other words, $v$ is symmetric with respect 
to the hyperplane $\{x_n=0\}$. The same analysis holds for any 
hyperplane containing $\vec{e_1}$. Coming back to the initial function 
$u$, this complete the proof of Theorem 1.2.

\section{Test-functions estimates}
We first introduce some definitions and notations. We consider a
family $(\ae)_{\eps>0}\in C^1(\Omega)$ and a function $a\in C^1(\Omega)$
such that there exists an open subset ${\mathcal U}\subset\rn$ such that
$\ae, a$ can be extended to ${\mathcal U}$ by $C^1-$functions that we
still denote by $\ae,a$. We assume that they satisfy
\bequa\label{hyp:ae}
\overline{\Omega}\subset\subset{\mathcal U}\hbox{ and }\lim_{\eps\to
0}\ae=a\hbox{ in }C^1_{loc}({\mathcal U}).\
\eequa
We assume that
\bequa\label{coerc:a}
\Delta+a\hbox{ is coercive in }\Omega,
\eequa
that is, there exists $c_0>0$ such that
$$\int_{\Omega}(|\nabla \varphi|^2+a\varphi^2)\, dx\geq
c_0\int_\Omega\varphi^2\, dx$$
for all $\varphi\in C^1_c(\Omega)$, the set of $C^1$-functions compactly
supported in $\Omega$. Finally, we let
$$\mu_{s,a}(\Omega)=\inf\frac{\int_{\Omega}(|\nabla u|^2+au^2)\,
dx}{\left(\int_{\Omega}\frac{|u|^{\crit}}{|x|^s}\,
dx\right)^{\frac{2}{\crit}}}.$$
Note that $\mu_{s,0}(\Omega)=\mu_{s}(\Omega)$. We let
$x_0\in\partial\Omega$. Since $\partial\Omega$ is smooth and $x_0\in
\partial\Omega$, there exist $U,V$ open subsets of $\rn$, there exists $I$
an open intervall of $\rr$, there exists $U'$ an open subset of
$\rr^{n-1}$ such that $0\in U=I\times U'$ and $x_0\in V$. There exist
$\varphi\in C^\infty(U,V)$ and $\varphi_0\in C^\infty(U')$ such that
\bequa\begin{array}{ll}
(i) & \varphi: U\to V\hbox{ is a }C^\infty-\hbox{diffeomorphism}\\
(ii) & \varphi(0)=x_0\\
(iii) & D_{0}\varphi=Id_{\rn}\\
(iv) & \varphi(U\cap\{x_1<0\})=\varphi(U)\cap\Omega\hbox{ and
}\varphi(U\cap\{x_1=0\})=\varphi(U)\cap\partial\Omega.\\
(v) & \varphi(x_1,y)=x_0+(x_1+\varphi_0(y),y)\hbox{ for all }(x_1,y)\in
I\times U'=U\\
(vi) & \varphi_0(0)=0\hbox{ and }\nabla\varphi_0(0)=0.
\end{array}\label{def:vphi}
\eequa
Here $D_x\varphi$ denotes the differential of $\varphi$ at $x$. This chart
will be useful throughout all the paper.

\medskip\noindent The first result we prove is an upper bound for
$\mu_{s,a}(\Omega)$.
\begin{prop}\label{prop:test}
Let $\Omega$ be a smooth bounded domain of $\rn$, $n\geq 3$, such that
$0\in\partial\Omega$. If $a\in C^0(\overline{\Omega})$ and $s\in
(0,2)$, then
$\mu_{s,a}(\Omega)\leq \mu_{s}(\rnm).$
 \end{prop}
\begin{proof} Let $\alpha>0$ and $u\in C_c^\infty(\rnm)\setminus\{0\}$
such that
$$\frac{\int_{\rnm}|\nabla u|^2\,
dx}{\left(\int_{\rnm}\frac{|u|^{\crit}}{|x|^s}\,
dx\right)^{\frac{2}{\crit}}}\leq \mu_{s}(\rnm)+\alpha.$$
Taking $x_0=0$ in (\ref{def:vphi}), we define
$$\ue(x)=\eps^{-\frac{n-2}{2}}u\left(\frac{\varphi^{-1}(x)}{\eps}\right)$$
for all $x\in\Omega$ and all $\eps>0$. As easily checked, for $\eps>0$
small enough, we have that
$$\ue\in C_c^\infty(\Omega).$$
With a change of variable, we get that
$$\int_\Omega\frac{|\ue|^{\crit}}{|x|^s}\,
dx=\int_{\rnm}\frac{\left|u(y)\right|^{\crit}}{\left|\frac{\varphi(\eps
y)}{\eps}\right|^s}\cdot |\hbox{Jac}(\eps y)|\, dy$$
Since $u$ is compactly supported, we get with point (iii) of
(\ref{def:vphi}) and Lebesgue's convergence theorem that
$$\lim_{\eps\to 0}\int_\Omega\frac{|\ue|^{\crit}}{|x|^s}\,
dx=\int_{\rnm}\frac{|u|^{\crit}}{|x|^s}\, dx.$$
On the other hand, we have that
$$\int_\Omega (|\nabla\ue|^2+a\ue^2)\, dx=\int_{\rnm}(|\nabla
u|_{g_\eps}^2+\eps^2 a\circ\varphi(\eps x) u^2)\cdot \sqrt{|g_\eps|}\,
dx,$$
where $(g_\eps(x))_{ij}=(\partial_i\varphi(\eps x),\partial_j\varphi(\eps
x))$, and $|g_\eps|=\det (g_\eps)$. With point (iii) of (\ref{def:vphi})
and Lebesgue's convergence theorem, we get that
$$\lim_{\eps\to 0}\int_\Omega (|\nabla\ue|^2+a\ue^2)\,
dx=\int_{\rnm}|\nabla u|^2\, dx.$$
As a consequence, we get that
$$\mu_{s,a}(\Omega)\leq \frac{\int_\Omega (|\nabla\ue|^2+a\ue^2)\,
dx}{\left(\int_\Omega\frac{|\ue|^{\crit}}{|x|^s}\,
dx\right)^{\frac{2}{\crit}}}=\frac{\int_{\rnm} |\nabla u|^2\,
dx}{\left(\int_{\rnm}\frac{|u|^{\crit}}{|x|^s}\,
dx\right)^{\frac{2}{\crit}}}+o(1)\leq \mu_{s}(\rnm)+\alpha+o(1)$$
where $\lim_{\eps\to 0}o(1)=0$. Letting $\eps\to 0$ and $\alpha\to 0$
yields the conclusion of the proposition.
\end{proof}

\section{The subcritical case}

\medskip\noindent{\bf Step 4.1:} In order to construct minimizers for
$\mu_{s,a}(\Omega)$, we consider a subcritical minimization problem for
which we recover compactness. This is the object of the following
proposition.
\begin{prop}\label{prop:subcrit} Let $\Omega$ be a smooth bounded domain
of $\rn$, $n\geq 3$ and $s\in (0,2)$. For any $\eps\in (0,\crit-2)$,
we let $\ae\in C^1(\overline{\Omega})$ such that $\Delta+\ae$ is coercive.
Then for any $\eps\in (0,\crit-2)$, the infimum
$$\mu_{s,\ae}^\eps(\Omega):=\inf_{u\in\huno\setminus\{0\}}\frac{\int_\Omega
(|\nabla u|^2+\ae u^2)\,
dx}{\left(\int_\Omega\frac{|u|^{\crit-\eps}}{|x|^s}\,
dx\right)^{\frac{2}{\crit-\eps}}},$$
is achieved by a function $\ue\in\huno$, $\ue> 0$ in $\Omega$. Moreover,
$\ue\in C^1(\overline{\Omega})\cap C^2(\overline{\Omega}\setminus\{0\})$
and can be assumed to satisfy the system
$$\left\{\begin{array}{ll}
\Delta\ue+\ae\ue=\frac{\ue^{\crit-1-\eps}}{|x|^s} & \hbox{ in }{\mathcal
D}'(\Omega)\\
\ue>0 & \hbox{ in }\Omega\\
\int_\Omega\frac{|u|^{\crit-\eps}}{|x|^s}\,
dx=(\mu_{s,\ae}^\eps(\Omega))^{\frac{\crit-\eps}{\crit-2-\eps}}&
\end{array}\right.$$
 \end{prop}

\begin{proof} This result is quite standard. We prove the proposition for
the sake of completeness. We claim that there exists a minimizer for
$\mu_{s,\ae}^\eps(\Omega)$. Indeed, let $(u_k)_{k\in\mathbb{N}}\in\huno$
be a minimizing sequence for $\mu_{s,\ae}^\eps(\Omega)$ such that
$$\int_{\Omega}\frac{|u_k|^{\crit-\eps}}{|x|^s}\, dx=1\hbox{ and
}\mu_{s,\ae}^\eps(\Omega)=\int_\Omega (|\nabla u_k|^2+\ae u_k^2)\,
dx+o(1)$$
where $\lim_{k\to +\infty}o(1)=0$. Since $\Vert u_k\Vert_{\huno}=O(1)$
when $k\to +\infty$, there exists $\tue\in\huno$ such that, up to a
subsequence, $u_k\rightharpoonup \tue$ weakly in $\huno$ when $k\to
+\infty$ and $\lim_{k\to +\infty}u_k(x)=\tue(x)$ a.e. in $\Omega$. Let
$\theta_k=u_k-\tue\in\huno$. As easily checked, we have that
\bequa\label{eq:lim:uk}
\mu_{s,\ae}^\eps(\Omega)=\int_\Omega(|\nabla \tue|^2+\ae\tue^2)\,
dx+\int_\Omega|\nabla \theta_k|^2\, dx+o(1),
\eequa
where $\lim_{k\to +\infty}o(1)=0$. Let $\eta\in C^\infty_c(\rr)$ such that
$\eta(x)=1$ for all $x\in [-1,1]$. Let $A>0$. With Lebesgue's theorem, we
have that
\beq
&&\left|\int_{\Omega}\frac{|u_k|^{\crit-\eps}}{|x|^s}\,
dx-\int_{\Omega}\frac{|\tue|^{\crit-\eps}}{|x|^s}\, dx\right|\\
&&=\left|\int_\Omega
\left(\eta\left(\frac{u_k}{A}\right)\frac{|u_k|^{\crit-\eps}}{|x|^s}-\eta\left(\frac{\tue}{A}\right)\frac{|\tue|^{\crit-\eps}}{|x|^s}\right)\,
dx\right|\\
&&+\int_\Omega
\left|1-\eta\left(\frac{u_k}{A}\right)\right|\frac{|u_k|^{\crit-\eps}}{|x|^s}\,
dx+\int_\Omega
\left|1-\eta\left(\frac{\tue}{A}\right)\right|\frac{|\tue|^{\crit-\eps}}{|x|^s}\,
dx\\
&&\leq  o(1)+\frac{1}{A^\eps}\int_\Omega
\left|1-\eta\left(\frac{u_k}{A}\right)\right|\frac{|u_k|^{\crit}}{|x|^s}\,
dx +\frac{1}{A^\eps}\int_\Omega
\left|1-\eta\left(\frac{\tue}{A}\right)\right|\frac{|\tue|^{\crit}}{|x|^s}\,
dx\\
&&\leq  o(1)+\frac{1}{A^\eps}\int_\Omega \frac{|u_k|^{\crit}}{|x|^s}\, dx
+\frac{1}{A^\eps}\int_\Omega \frac{|\tue|^{\crit}}{|x|^s}\, dx\\
&&\leq  o(1)+\frac{1}{A^\eps}\mu_{s}(\Omega)^{-\frac{\crit}{2}}\left(\Vert
u_k\Vert_{\huno}^{\crit}+\Vert \tue\Vert_{\huno}^{\crit}\right)
\eeq
where $\lim_{k\to +\infty}o(1)=0$. Letting $k\to +\infty$, and then $A\to
+\infty$, we get that
$$\lim_{k\to +\infty}\int_{\Omega}\frac{|u_k|^{\crit-\eps}}{|x|^s}\,
dx=\int_{\Omega}\frac{|\tue|^{\crit-\eps}}{|x|^s}\, dx.$$
It then follows that $\int_{\Omega}\frac{|\tue|^{\crit-\eps}}{|x|^s}\,
dx=1$. With the definition of $\mu_{s,\ae}^\eps(\Omega)$, we then get that
$$\mu_{s,\ae}^\eps(\Omega)\leq\int_\Omega(|\nabla \tue|^2+\ae\tue^2)\,
dx.$$
With (\ref{eq:lim:uk}), we then get that $\lim_{k\to +\infty}\theta_k=0$
in $\huno$. As a consequence, $\mu_{s,\ae}^\eps(\Omega)$ is attained by
$\tue$. This proves the claim.

\medskip\noindent Up to replacing $\tue$ by $|\tue|$, we can assume that
$\tue \geq 0$. We let
$$\ue=\mu_{s,\ae}^\eps(\Omega)^{\frac{1}{\crit-2-\eps}}\tue.$$
As easily checked, $\ue\geq 0$ is also a minimizer for
$\mu_{s,\ae}^\eps(\Omega)$. It satisfies
$$\Delta\ue+\ae\ue=\frac{\ue^{\crit-1-\eps}}{|x|^s} \hbox{ in }{\mathcal
D}'(\Omega).$$
Moreover, it follows from the appendix and standard elliptic theory that
$\ue\in C^1(\overline{\Omega})\cap C^2(\overline{\Omega}\setminus\{0\})$.
Since $\Delta\ue\geq 0$ in $\Omega$ and $\ue\not\equiv 0$, it follows from
the strong comparison principle that $\ue>0$ in $\Omega$.
\end{proof}

\medskip\noindent{\bf Step 4.2: } For any $\eps\in (0,\crit-2)$, we let
$(\ae), a$ as in (\ref{hyp:ae}) and (\ref{coerc:a}). We let
$\mu_{s,\ae}^\eps(\Omega)$ as in Proposition \ref{prop:subcrit}. We claim
that
$$\lim_{\eps\to 0}\mu_{s,\ae}^\eps(\Omega)=\mu_{s,a}(\Omega).$$
Indeed, we let $\alpha>0$ and let $u\in C^\infty_c(\Omega)\setminus\{0\}$
such that
$$\frac{\int_\Omega (|\nabla u|^2+au^2)\,
dx}{\left(\int_\Omega\frac{|u|^{\crit}}{|x|^s}\,
dx\right)^{\frac{2}{\crit}}}\leq \mu_{s,a}(\Omega)+\alpha.$$
We have that
$$\lim_{\eps\to 0}\frac{\int_\Omega (|\nabla u|^2+\ae u^2)\,
dx}{\left(\int_\Omega\frac{|u|^{\crit-\eps}}{|x|^s}\,
dx\right)^{\frac{2}{\crit-\eps}}}=\frac{\int_\Omega (|\nabla u|^2+au^2)\,
dx}{\left(\int_\Omega\frac{|u|^{\crit}}{|x|^s}\,
dx\right)^{\frac{2}{\crit}}}\leq \mu_{s,a}(\Omega)+\alpha.$$
Letting $\eps\to 0$ and $\alpha\to 0$, we get that
\bequa\label{ineq:cv:mue:1}
\limsup_{\eps\to 0}\mu_{s,\ae}^\eps(\Omega)\leq \mu_{s,a}(\Omega).
\eequa
We now let $v\in C^\infty_c(\Omega)\setminus\{0\}$. It follows from
H\"older's inequality that
$$\left(\int_\Omega\frac{|v|^{\crit-\eps}}{|x|^s}\,
dx\right)^{\frac{2}{\crit-\eps}}\leq
\left(\int_\Omega\frac{dx}{|x|^s}\right)^{\frac{2\eps}{\crit\cdot(\crit-\eps)}}\left(\int_\Omega\frac{|v|^{\crit}}{|x|^s}\,
dx\right)^{\frac{2}{\crit}}$$
and then
\beq
\frac{\int_\Omega (|\nabla v|^2+a v^2)\,
dx}{\left(\int_\Omega\frac{|v|^{\crit}}{|x|^s}\,
dx\right)^{\frac{2}{\crit}}}&\leq&
\left(\int_\Omega\frac{dx}{|x|^s}\right)^{\frac{2\eps}{\crit\cdot(\crit-\eps)}}\cdot
\frac{\int_\Omega (|\nabla v|^2+\ae v^2)\,
dx}{\left(\int_\Omega\frac{|v|^{\crit-\eps}}{|x|^s}\,
dx\right)^{\frac{2}{\crit-\eps}}}\\
&&+\frac{\int_\Omega (a-\ae)v^2\,
dx}{\left(\int_\Omega\frac{|v|^{\crit-\eps}}{|x|^s}\,
dx\right)^{\frac{2}{\crit-\eps}}}
\eeq
for $\eps>0$ small. Here, we have used that $\Delta+\ae$ is coercive on
$\Omega$ for $\eps>0$ small, which is a consequence of (\ref{hyp:ae}) and
(\ref{coerc:a}). Taking the infimum, using H\"older's inequality and
(\ref{hyp:ae}), we get that
\bequa\label{ineq:cv:mue:2}
\mu_{s,a}(\Omega)\leq (1+o(1))\mu_{s,\ae}^\eps(\Omega)
\eequa
where $\lim_{\eps\to 0}o(1)=0$. The conclusion of  Step 4.2 then follows
from (\ref{ineq:cv:mue:1}) and (\ref{ineq:cv:mue:2}).

\medskip\noindent{\bf  Step 4.3:} We prove that, when it is nonzero, the
weak limit of the $\ue$'s is a minimizer for $\mu_{s,a}(\Omega)$. This is
the object of the following proposition.
\begin{prop}\label{prop:min:nonzero}
Let $\Omega$ be a smooth bounded domain of $\rn$, $n\geq 3$, such that
$0\in\partial\Omega$. For $s\in (0,2)$ and $\eps\in (0,\crit-2)$,
we let $\ae, a$ be as in (\ref{hyp:ae}) and (\ref{coerc:a}). For any 
$\eps\in
(0,\crit-2)$, let $\mu_{s,\ae}^\eps(\Omega)$ and $\ue$ be as in 
Proposition
\ref{prop:subcrit}. Then there exists $u_0\in \huno$ such that, up to a
subsequence, $\ue\rightharpoonup u_0$ weakly in $\huno$ when $\eps\to 0$.
If $u_0\not\equiv 0$, then $\lim_{\eps\to 0}\ue=u_0$ strongly in $\huno$
and $u_0$ is a minimizer for $\mu_{s,a}(\Omega)$. In particular,
$\mu_{s,a}(\Omega)$ is attained.
\end{prop}

\begin{proof} It is clear from Proposition \ref{prop:subcrit} and the
hypothesis (\ref{hyp:ae}) and (\ref{coerc:a}) that
$$\Vert\ue\Vert_{\huno}=O(1)$$
when $\eps\to 0$. Then there exists $u_0\in \huno$ such that, up to a
subsequence, $\ue\rightharpoonup u_0$ weakly in $\huno$ when $\eps\to 0$.
We assume that $u_0\not\equiv 0$. It then follows from the definition of
$\mu_{s,a}(\Omega)$ that
$$\frac{\int_\Omega (|\nabla u_0|^2+a u_0^2)\,
dx}{\left(\int_\Omega\frac{|u_0|^{\crit}}{|x|^s}\,
dx\right)^{\frac{2}{\crit}}}\geq \mu_{s,a}(\Omega).$$
Testing the weak inequality
$\Delta\ue+\ae\ue=\frac{\ue^{\crit-1-\eps}}{|x|^s}$ on $u_0$ and letting
$\eps\to 0$, we get that
$$\int_\Omega (|\nabla u_0|^2+a u_0^2)\,
dx=\int_\Omega\frac{|u_0|^{\crit}}{|x|^s}\, dx.$$
We then obtain that
$$\int_\Omega\frac{|u_0|^{\crit}}{|x|^s}\, dx\geq
\mu_{s,a}(\Omega)^{\frac{\crit}{\crit-2}}.$$
Since $\ue\rightharpoonup u_0$ when $\eps\to 0$, we get with the
definition of $\ue$ in Proposition \ref{prop:subcrit} and  Step 4.2 that
$$\int_\Omega\frac{|u_0|^{\crit}}{|x|^s}\, dx\leq \liminf_{\eps\to
0}\int_\Omega\frac{|\ue|^{\crit-\eps}}{|x|^s}\,
dx=\mu_{s,a}(\Omega)^{\frac{\crit}{\crit-2}}.
$$
Consequently, we get that
\bequa\label{eq:lim:uzero:1}
\int_\Omega (|\nabla u_0|^2+a u_0^2)\,
dx=\int_\Omega\frac{|u_0|^{\crit}}{|x|^s}\,
dx=\mu_{s,a}(\Omega)^{\frac{\crit}{\crit-2}}.
\eequa
Since
$\mu_{s,\ae}^\eps(\Omega)^{\frac{\crit-\eps}{\crit-2-\eps}}=\int_\Omega
(|\nabla \ue|^2+\ae\ue^2)\, dx$, we get with the definition of $\ue$ in
Proposition \ref{prop:subcrit} that
\bequa\label{eq:lim:uzero:2}
\mu_{s,a}(\Omega)^{\frac{\crit}{\crit-2}}=\int_\Omega (|\nabla
u_0|^2+au_0^2)\, dx+\int_\Omega |\nabla (\ue-u_0)|^2\, dx+o(1)
\eequa
with $\lim_{\eps\to 0}o(1)=0$. It follows from (\ref{eq:lim:uzero:1}) and
(\ref{eq:lim:uzero:2}) that $\lim_{\eps\to 0}\ue=u_0$ in $\huno$. As
easily checked, in this case, $u_0$ is a minimizer for
$\mu_{s,a}(\Omega)$.
\end{proof}

\section{Preliminary Blow-Up analysis}
From now on, we let $\Omega$ be a smooth bounded domain of $\rn$, $n\geq
3$, such that $0\in\partial\Omega$. We let $s\in (0,2)$. For any $\eps>0$,
we let $\pe\in [0,\crit-2)$ such that
\bequa\label{lim:pe}
\lim_{\eps\to 0}\pe=0.
\eequa
We let $a\in C^1(\overline{\Omega})$ and a family $(\ae)_{\eps>0}\in
C^1(\overline{\Omega})$ such that (\ref{hyp:ae}) and (\ref{coerc:a}) hold.
For any $\eps>0$, we consider $\ue\in \huno\cap
C^2(\overline{\Omega}\setminus\{0\})$ a solution to the system
\bequa\label{syst:ue}
\left\{\begin{array}{ll}
\Delta\ue+\ae\ue=\frac{\ue^{\crit-1-\pe}}{|x|^s}& \hbox{ in }{\mathcal
D}'(\Omega)\\
\ue>0&\hbox{ in }\Omega
\end{array}\right.
\eequa
for all $\eps>0$. We assume that $\ue$ is of minimal energy type, that is
\bequa\label{hyp:nrj:min}
\int_\Omega\frac{|\ue|^{\crit-\pe}}{|x|^s}\,
dx=\mu_{s}(\Omega)^{\frac{\crit}{\crit-2}}+o(1)
\eequa
where $\lim_{\eps\to 0}o(1)=0$. Note that it follows from (\ref{hyp:ae}),
(\ref{coerc:a}), (\ref{syst:ue}) and (\ref{hyp:nrj:min}) that
\bequa\label{bnd:ue}
\Vert\ue\Vert_{\huno}=O(1)
\eequa
when $\eps\to 0$. We also assume that blow-up occurs, that is
\bequa\label{hyp:blowup}
\ue\rightharpoonup 0
\eequa
weakly in $\huno$ when $\eps\to 0$. Such a family arises naturally when
$u_0\equiv 0$ in Propositions \ref{prop:subcrit} and
\ref{prop:min:nonzero}. In the remaining sections, we describe precisely
the behaviour of the $\ue$'s. We follow the strategy developed in
\cite{dhr}.

\medskip\noindent It follows from Proposition \ref{prop:app} of the
Appendix that $\ue\in C^0(\overline{\Omega})$. We let $\xe\in\Omega$ and
$\me,\ke>0$ such that
\bequa\label{def:me:xe}
\max_\Omega\ue=\ue(\xe)=\me^{-\frac{n-2}{2}}\hbox{ and
}\ke:=\me^{1-\frac{\pe}{\crit-2}}.
\eequa
We let $\varphi:U\to V$ a local chart as in (\ref{def:vphi}) with $x_0=0$,
where $U,V$ are open neighborhoods of $0$. For any $\eps>0$ and any $x\in
\frac{U}{\ke}\cap\{x_1\leq 0\}$, we define the maximum rescaling of $\ue$
as follows
\bequa\label{def:resc:ue}
\ve(x):=\frac{\ue\circ\varphi(\ke x)}{\ue(\xe)},
\eequa
where $\xe,\ke$ are as in (\ref{def:me:xe}). As easily checked, for any
$\eta\in C_c^\infty(\rn)$, we have that $\eta\ve\in\hunrnm$. In this
section, we prove the following proposition:
\begin{prop}\label{prop:sec3:3}
Let $\Omega$ be a smooth bounded domain of $\rn$, $n\geq 3$ and $s\in
(0,2)$. Consider $(\pe)_{\eps>0}$ such that $\pe\in [0,\crit-2)$ for all
$\eps>0$. We consider $(\ue)_{\eps>0}\in \huno$ such that (\ref{hyp:ae}),
(\ref{coerc:a}), (\ref{syst:ue}), (\ref{hyp:nrj:min}) and
(\ref{hyp:blowup}) hold.  Let $\ve$ be as in (\ref{def:resc:ue}). Then
there exists $v\in \hunrnm\setminus\{0\}$ such that for any $\eta\in
C_c^\infty(\rn)$,
$$\eta\ve\rightharpoonup\eta v\hbox{ in }\hunrnm$$
when $\eps\to 0$. Moreover, $v$ verifies that
$$\Delta v=\frac{v^{\crit-1}}{|x|^s}\hbox{ in }{\mathcal D}'(\rnm)$$
and
$$\int_{\rnm}|\nabla v|^2\,
dx=\mu_{s}(\Omega)^{\frac{\crit}{\crit-2}}=\mu_{s}(\rnm)^{\frac{\crit}{\crit-2}}.$$
In addition, there exists $\theta\in (0,1)$ such that $v\in
C^{1,\theta}(\overline{\rnm})$ and
$$\ve\to v\hbox{ in }C^{1,\theta}_{loc}(\overline{\rnm})$$
when $\eps\to 0$. Moreover, we have that
\bequa\label{lim:re}
\lim_{\eps\to 0}\me^{\pe}=1.
\eequa
\end{prop}
\begin{proof} Steps 5.1 to 5.9 below are devoted to the proof of this
Proposition.

\medskip\noindent{\bf  Step 5.1:} We claim that
\bequa\label{concl:step41}
\me=o(1)
\eequa
when $\eps\to 0$. We proceed by contradiction and assume that
$\lim_{\eps\to 0}\me\neq 0$. In this case, up to a subsequence, there
exists $C>0$ such that $\ue(x)\leq C$ for all $x\in\Omega$ and all
$\eps>0$. Since (\ref{hyp:blowup}) hold, it follows from standard elliptic
theory (see for instance \cite{gt}) that $\lim_{\eps\to 0}\ue=0$ in
$C^0(\overline{\Omega})$. A contradiction with (\ref{hyp:nrj:min}). This
proves (\ref{concl:step41}).

\medskip\noindent{\bf  Step 5.2:} We claim that
\bequa\label{bnd:xe:ke}
|\xe|=O(\ke)
\eequa
when $\eps\to 0$. We proceed by contradiction and assume that
\bequa\label{lim:step:1}
\lim_{\eps\to 0}\frac{|\xe|}{\ke}=+\infty.
\eequa
For any $\eps>0$, we let
\bequa\label{def:be}
\be=|\xe|^{\frac{s}{2}}\ue(\xe)^{\frac{2+\pe-\crit}{2}}=|\xe|^{\frac{s}{2}}\ke^{\frac{2-s}{2}}.
\eequa
It follows from the definition (\ref{def:be}) of $\be$ and
(\ref{lim:step:1}) that
\bequa\label{ppty:be}
\lim_{\eps\to 0}\be=0,\; \lim_{\eps\to 0}\frac{\be}{\ke}=+\infty\hbox{ and
}\lim_{\eps\to 0}\frac{\be}{|\xe|}=0
\eequa
when $\eps\to 0$.

\medskip\noindent{\it Case 5.2.1:} We assume that there exists $\rho>0$
such that
$$\frac{d(\xe,\partial\Omega)}{\be}\geq 2\rho$$
for all $\eps>0$. For $x\in B_{2\rho}(0)$ and $\eps>0$, we define
$$\bve(x):=\frac{\ue(\xe+\be x)}{\ue(\xe)}.$$
Note that this is well defined since $\xe+\be x\in\Omega$ for all $x\in
B_{2\rho}(0)$. As easily checked, we have that
$$\Delta\bve+\be^2\ae(\xe+\be
x)\bve=\frac{\bve^{\crit-1-\pe}}{\left|\frac{\xe}{|\xe|}+\frac{\be}{|\xe|}\cdot
x\right|^s}$$
weakly in $B_{2\rho}(0)$. Since (\ref{ppty:be}) holds, we have that
$$\Delta\bve+\be^2\ae(\xe+\be x)\bve=(1+o(1))\bve^{\crit-1-\pe}$$
weakly in $B_{2\rho}(0)$, where $\lim_{\eps\to 0}o(1)=0$ in
$C^0_{loc}(B_{2\rho}(0))$. Since $0\leq \bve(x)\leq \bve(0)=1$ for all
$x\in B_{2\rho}(0)$, it follows from standard elliptic theory that there
exists $v\in C^1(B_{2\rho}(0))$ such that $v\geq 0$ and
$$\bve\to \overline{v}$$
in $C^1_{loc}(B_{2\rho}(0))$ when $\eps\to 0$. In particular,
\bequa\label{lim:v:case1:nonzero}
\overline{v}(0)=\lim_{\eps\to 0}\bve(0)=1.
\eequa
With a change of variables and the definition (\ref{def:be}) of $\be$, we
get that
\beq
&&\int_{\Omega\cap B_{\rho\be}(\xe)}\frac{\ue^{\crit-\pe}}{|x|^s}\,
dx=\frac{\ue(\xe)^{\crit-\pe}\be^n}{|\xe|^s}\int_{B_{\rho}(0)}\frac{\bve^{\crit-\pe}}{\left|\frac{\xe}{|\xe|}+\frac{\be}{|\xe|}\cdot
x\right|^s}\, dx\\
&&\geq
\left(\frac{\be}{\ke}\right)^{n-2}\int_{B_{\rho}(0)}\frac{\bve^{\crit-\pe}}{\left|\frac{\xe}{|\xe|}+\frac{\be}{|\xe|}x\right|^s}\,
dx.
\eeq
Using (\ref{hyp:nrj:min}), (\ref{ppty:be}) and passing to the limit
$\eps\to 0$ (note that $\me^{-1}\geq 1$ for $\eps>0$ small), we get that
$$\int_{B_{\rho}(0)}\overline{v}^{\crit}\, dx=0,$$
and then $\overline{v}\equiv 0$ in $B_\rho(0)$. A contradiction with
(\ref{lim:v:case1:nonzero}). Then (\ref{lim:step:1}) does not hold. This
proves that (\ref{bnd:xe:ke}) holds in Case 5.2.1.

\medskip\noindent{\it Case 5.2.2:} We assume that, up to a subsequence,
\bequa\label{lim:d:be:0}
\lim_{\eps\to 0}\frac{d(\xe,\partial\Omega)}{\be}=0.
\eequa
In this case,
$$\lim_{\eps\to 0}\xe=x_0\in\partial\Omega.$$
Since $x_0\in\partial\Omega$, we let $\varphi:U\to V$ as in
(\ref{def:vphi}), where $U,V$ are open neighborhoods of $0$ and $x_0$
respectively. We let $\tue=\ue\circ\varphi$, which is defined on $U\cap
\{x_1\leq 0\}$. For any $i,j=1,...,n$, we let
$g_{ij}=(\partial_i\varphi,\partial_j\varphi)$, where $(\cdot,\cdot)$
denotes the Euclidean scalar product on $\rn$, and we consider $g$ as a
metric on $\rn$. We let $\Delta_g=-div_g(\nabla)$ the Laplace-Beltrami
operator with respect to the metric $g$. In our basis, we have that
$$\Delta_g=-g^{ij}\left(\partial_{ij}-\Gamma_{ij}^k\partial_k\right),$$
where $g^{ij}=(g^{-1})_{ij}$ are the coordinates of the inverse of the
tensor $g$ and the $\Gamma_{ij}^k$'s are the Christoffel symbols of the
metric $g$. As easily checked, we have that
$$\Delta_g\tue+\ae\circ\varphi(x)\cdot
\tue=\frac{\tue^{\crit-1-\pe}}{|\varphi(x)|^s}$$
weakly in $U\cap\{x_1<0\}$. We let $\ze\in\partial\Omega$ such that
\bequa\label{def:ze}
|\ze-\xe|=d(\xe,\partial\Omega).
\eequa
We let $\txe,\tze\in U$ such that
\bequa\label{def:txe:tze}
\varphi(\txe)=\xe\hbox{ and }\varphi(\tze)=\ze.
\eequa
It follows from the properties (\ref{def:vphi}) of $\varphi$ that
\bequa\label{ppty:txe:tze}
\lim_{\eps\to 0}\txe=\lim_{\eps\to 0}\tze=0,\; (\txe)_1<0\hbox{ and
}(\tze)_1=0.\eequa
At last, we let
$$\tve(x):=\frac{\tue(\tze+\be x)}{\tue(\txe)}$$
for all $x\in \frac{U-\tze}{\be}\cap \{x_1<0\}$. With
(\ref{ppty:txe:tze}), we get that $\tve$ is defined on
$B_R(0)\cap\{x_1<0\}$ for all $R>0$, as soon as $\eps$ is small enough.
The function $\tve$ verifies
$$\Delta_{\tge}\tve+\be^2\ae\circ\varphi(\tze+\be
x)\tve=\frac{\tve^{\crit-1-\pe}}{\left|\frac{\varphi(\tze+\be
x)}{|\xe|}\right|^s}$$
weakly in $B_R(0)\cap\{x_1<0\}$. In this expression, $\tge=g(\tze+\be x)$
and $\Delta_{\tge}$ is the Laplace-Beltrami operator with respect to the
metric $\tge$. With (\ref{lim:d:be:0}), (\ref{def:ze}) and
(\ref{def:txe:tze}), we get that
$$\varphi(\tze+\be x)=\xe+O_R(1)\be,$$
for all $x\in B_{R}(0)\cap\{x_1\leq 0\}$ and all $\eps>0$, where there
exists $C_R>0$ such that $|O_R(1)|\leq C_R$ for all $x\in
B_{R}(0)\cap\{x_1<0\}$. With (\ref{ppty:be}), we then get that
$$\lim_{\eps\to 0}\frac{|\varphi(\tze+\be x)|}{|\xe|}=1$$
in $C^0(B_{R}(0)\cap\{x_1\leq 0\})$. It then follows that
$$\Delta_{\tge}\tve+\be^2\ae\circ\varphi(\tze+\be
x)\tve=(1+o(1))\tve^{\crit-1-\pe}$$
weakly in $B_R(0)\cap\{x_1<0\}$, where $\lim_{\eps\to 0}o(1)=0$ in
$C^0(B_R(0)\cap\{x_1\leq 0\})$. Since $\tve$ vanishes on
$B_R(0)\cap\{x_1=0\}$ (in the sense of the trace) and that $0\leq \tve\leq
1$, it follows from standard elliptic theory that there exists $\tv\in
C^1(B_R(0)\cap\{x_1\leq 0\})$ such that $$\lim_{\eps\to 0}\tve=\tv$$
in $C^{0}(B_{\frac{R}{2}}(0)\cap \{x_1\leq 0\})$. In particular,
\bequa\label{eq:v:vanish}
\tv\equiv 0\hbox{ on }B_{\frac{R}{2}}(0)\cap \{x_1=0\}.
\eequa
Moreover, it follows from (\ref{def:ze}) and (\ref{def:txe:tze}) that
$$\tve\left(\frac{\txe-\tze}{\be}\right)=1\hbox{ and }\lim_{\eps\to
0}\frac{\txe-\tze}{\be}=0.$$
In particular, $\tv(0)=1$. A contradiction with (\ref{eq:v:vanish}). Then
(\ref{lim:step:1}) does not hold. This proves (\ref{bnd:xe:ke}) in Case
5.2.2.

\smallskip\noindent In both cases, we have obtained that (\ref{bnd:xe:ke})
holds. This concludes  Step 5.2.

\medskip\noindent A consequence of (\ref{bnd:xe:ke}) is that
$\lim_{\eps\to 0}\xe=0\in\partial\Omega$. We let $\varphi:U\to V$ as in
(\ref{def:vphi}) be a local chart of $\partial\Omega$ with $x_0=0$ (in
other words, $\varphi(0)=0$), where $U,V$ are open neighborhoods of $0$.
We write
$$\xe=\varphi(\xune, \ze),$$
where $\xune<0$ and $\ze\in \rr^{n-1}$ are such that $(\xune,\ze)\in U$.

\medskip\noindent{\bf  Step 5.3:} We claim that
\bequa\label{lim:d:ze}
d(\xe,\partial\Omega)=(1+o(1))|\xune|=O(\ke)\hbox{ and }\ze=O(\ke),
\eequa
when $\eps\to 0$. Indeed, with (\ref{bnd:xe:ke}), we get that
\bequa\label{lim:prop32:1}
d(\xe,\partial\Omega)\leq |\xe|=O(\ke)
\eequa
when $\eps\to 0$. We first remark that
$$d(\xe,\partial\Omega)\leq d(\xe,\varphi(0,\ze))=|\xune|.$$
We let $\ae\in \hbox{span}(\vec{e}_2,...,\vec{e}_n)$ and
$Y_\eps=\varphi(0,a_\eps)\in\partial\Omega$ such that
$d(\xe,\partial\Omega)=|\xe-Y_\eps|$. Since $d(\xe,\partial\Omega)\leq
|\xune|$, we get that
$$\ze-a_\eps=O(|\xune|),$$
when $\eps\to 0$. Since $\nabla\varphi_0(0)=0$ (where $\varphi_0$ is as in
(\ref{def:vphi})), we get that
$$\varphi_0(\ze)=\varphi_0(a_\eps)+o(|\ze-a_\eps|)=\varphi_0(a_\eps)+o(|\xune|)$$
when $\eps\to 0$. Moreover,
\beq
d(\xe,\partial\Omega)&=&|\xe-Y_\eps|\\
&=& |(\xune+\varphi_0(\ze)-\varphi_0(a_\eps), \ze-a_\eps)|\\
&=& |(\xune+o(|\xune|), \ze-a_\eps)|\leq |\xune|
\eeq
when $\eps\to 0$. It then follows that $\ze-a_\eps=o(|\xune|)$ and
$d(\xe,\partial\Omega)=(1+o(1))|\xune|$ when $\eps\to 0$. This last
result, (\ref{bnd:xe:ke}) and (\ref{lim:prop32:1}) prove (\ref{lim:d:ze}).

\medskip\noindent{\bf  Step 5.4:} We let
\bequa\label{def:lae}
\lae:=-\frac{\xune}{\ke}>0\hbox{ and }\theta_\eps:=\frac{\ze}{\ke}.
\eequa
It follows from (\ref{lim:d:ze}) that there exist $\lambda_0\geq 0$ and
$\theta_0\in\rr^{n-1}$ such that
\bequa\label{lim:lae}
\lim_{\eps\to 0}\lae=\lambda_0\hbox{ and }\lim_{\eps\to
0}\theta_\eps=\theta_0.
\eequa
For any $\eps>0$ and any $x\in \frac{U}{\ke}\cap\{x_1\leq 0\}$, we let (as
in (\ref{def:resc:ue}))
\bequa\label{def:ve:sec3}
\ve(x):=\frac{\ue\circ\varphi(\ke x)}{\ue(\xe)},
\eequa
where $\varphi:U\to V$ is defined in (\ref{def:vphi}) (with $x_0=0$) and
$\ke,\xe$ are as in (\ref{def:me:xe}). As easily checked, for any $\eta\in
C_c^\infty(\rn)$, we have that
$$\eta\ve\in\hunrnm$$
for all $\eps>0$. We go on in the proof of Proposition \ref{prop:sec3:3}.

\medskip\noindent{\bf  Step 5.5:} We claim that for any $\eta\in
C_c^\infty(\rn)$, there exists $v_\eta\in \hunrnm$ such that, up to a
subsequence,
$$\eta\ve\rightharpoonup v_\eta$$
weakly in $\hunrnm$. Indeed, as easily checked, we have that
$$\nabla(\eta\ve)(x)=\ve\nabla\eta+\frac{\ke}{\ue(\xe)}\eta\cdot D_{(\ke
x)}\varphi[(\nabla\ue)(\varphi(\ke x))],$$
for all $\eps>0$ and all $x\in\rnm$. In this expression, $D_x\varphi$ is
the differential of the function $\varphi$ at $x$. It is standard that for
any $\alpha>0$, there exists $C_\alpha>0$ such that
$$(x+y)^2\leq C_\alpha x^2+(1+\alpha)\cdot y^2$$
for all $x,y>0$. With this inequality, we get that
\beq
&&\int_{\rnm}|\nabla(\eta\ve)|^2\, dx\leq
C_\alpha\int_{\rnm}|\nabla\eta|^2\ve^2\, dx\\
&&+(1+\alpha)\int_{\rnm}\eta^2\frac{\ke^2}{\ue(\xe)^2}\cdot|D_{(\ke
x)}\varphi[(\nabla\ue)(\varphi(\ke x))]|^2\, dx.
\eeq
Since $D_0\varphi=Id_{\rn}$, we get that with H\"older's inequality and a
change of variables that
\beqn
&&\int_{\rnm} |\nabla(\eta\ve)|^2\, dx\leq
C_\alpha\int_{\rnm}|\nabla\eta|^2\ve^2\, dx\nonumber\\
&&+(1+\alpha)\cdot(1+O(\ke))\int_{\rnm}\eta^2\frac{\ke^2}{\ue(\xe)^2}\cdot|\nabla\ue|^2(\varphi(\ke
x))\, dx\nonumber\\
&&\leq
C_\alpha\Vert\nabla\eta\Vert_{n}^2\cdot\Vert\ve\Vert_{L^{\frac{2n}{n-2}}(\hbox{Supp
}\nabla\eta)}^{2}\nonumber\\
&&+(1+\alpha)\cdot(1+O(\ke))\cdot\me^{\frac{\pe(n-2)}{\crit-2}}\int_\Omega|\nabla\ue|^2\,
dx\label{ineq:ve:1}
\eeqn
With another change of variables, we get that
\beqn
&&\int_{\rnm} |\nabla(\eta\ve)|^2\, dx\leq
C_\alpha\cdot\me^{\frac{(n-2)\pe}{\crit-2}}\Vert\nabla\eta\Vert_{n}^2\cdot\Vert\ue\Vert_{L^{\frac{2n}{n-2}}(\Omega)}^{2}\nonumber\\
&&+(1+\alpha)\cdot(1+O(\ke))\cdot\me^{\frac{\pe(n-2)}{\crit-2}}\int_\Omega|\nabla\ue|^2\,
dx.\label{ineq:ve:2}
\eeqn
With (\ref{bnd:ue}), Sobolev's inequality and since $\me^{\pe}\leq 1$ for
all $\eps>0$ small enough, we get with (\ref{ineq:ve:2}) that
$$\Vert\eta\ve\Vert_{\hunrnm}=O(1)$$
when $\eps\to 0$. It then follows that there exists $v_\eta\in\hunrnm$
such that, up to a subsequence, $\eta\ve\rightharpoonup v_\eta$ weakly in
$\hunrnm$ when $\eps\to 0$. This concludes  Step 5.5.

\medskip\noindent{\bf  Step 5.6:} We claim that there exists $v\in 
\hunrnm$
such that for any $\eta\in C_c^\infty(\rn)$, we have that, up to a
subsequence,
$$\eta\ve\rightharpoonup \eta v$$
weakly in $\hunrnm$ when $\eps\to 0$. Indeed, we let $\eta_1\in
C_c^\infty(\rn)$ such that $\eta_1\equiv 1$ in $B_1(0)$ and $\eta_1\equiv
0$ in $\rn\setminus B_2(0)$. For any $R>0$, we let
$\eta_R(x)=\eta_1(\frac{x}{R})$ for all $x\in\rn$. With a diagonal
argument, we can assume that, up to a subsequence, for any $R>0$, there
exists $v_R\in\hunrnm$ such that
$$\eta_R\ve\rightharpoonup v_R$$
weakly in $\hunrnm$ when $\eps\to 0$, and that $(\eta_R\ve)(x)\to v_R(x)$
when $\eps\to 0$ for a.e. $x\in\rnm$. Letting $\eps\to 0$ in
(\ref{ineq:ve:2}), with (\ref{bnd:ue}), Sobolev's inequality and since
$\me^{\pe}\leq 1$ for all $\eps>0$ small enough, we get that there exists
a constant $C>0$ independant of $R$ such that
$$\int_{\rnm} |\nabla v_R|^2\, dx\leq
C_\alpha\Vert\nabla\eta_R\Vert_{n}^2\cdot C+(1+\alpha)\cdot C$$
for all $R>0$. Since
$\Vert\nabla\eta_R\Vert_{n}^2=\Vert\nabla\eta_1\Vert_{n}^2$ for all $R>0$,
we get that there exists $C>0$ independant of $R$ such that
$$\int_{\rnm} |\nabla v_R|^2\, dx\leq C$$
for all $R>0$. It then follows that there exists $v\in \hunrnm$ such that
$v_R\rightharpoonup v$ weakly in $\hunrnm$ when $R\to +\infty$ and
$v_R(x)\to v(x)$ when $R\to +\infty$ for a.e. $x\in\rnm$. As easily
checked, we then obtain that $v_\eta=\eta v$ (we omit the proof of this
fact. It is straightforward). This ends  Step 5.6.

\medskip\noindent{\bf  Step 5.7:} We claim that
$$v\not\equiv 0.$$
Indeed, we let $R>0$. We proceed as in Case 5.2.2 of the proof of
(\ref{bnd:xe:ke}) in  Step 5.2, for any $i,j=1,...,n$, we let
$(\tge)_{ij}=(\partial_i\varphi(\ke x),\partial_j\varphi(\ke x))$, where
$(\cdot,\cdot)$ denotes the Euclidean scalar product on $\rn$. We consider
$\tge$ as a metric on $\rn$. We let
$$\Delta_{\tge}=-{\tge}^{ij}\left(\partial_{ij}-\Gamma_{ij}^k(\tge)\partial_k\right),$$
where $\tge^{ij}:=(\tge^{-1})_{ij}$ are the coordinates of the inverse of
the tensor $\tge$ and the $\Gamma_{ij}^k(\tge)$'s are the Christoffel
symbols of the metric $\tge$. With a change of variable and the definition
(\ref{def:ve:sec3}), equation (\ref{syst:ue}) rewrites as
\bequa\label{eq:ve:sec3}
\Delta_{\tge}(\eta_R\ve)+\ke^2\ae\circ\varphi(\ke
x)\eta_R\ve=\frac{(\eta_R\ve)^{\crit-1-\pe}}{\left|\frac{\varphi(\ke
x)}{\ke}\right|^s}\hbox{ in }{\mathcal D}'(B_R(0)\cap\{x_1<0\})
\eequa
for all $\eps>0$. With (\ref{def:me:xe}), (\ref{def:ve:sec3}) and since
$s\in (0,2)$, we get that $0\leq\ve\leq 1$ and that there exists
$p>\frac{n}{2}$ such that the RHS of (\ref{eq:ve:sec3}) is bounded in
$L^p$ when $\eps\to 0$. It then follows from standard elliptic theory that
there exists $\alpha>0$ such that
$$\Vert \eta_R\ve\Vert_{C^{0,\alpha}(B_{R/2}(0)\cap\{x_1\leq 0\})}=O(1)$$
when $\eps\to 0$. It then follows from Ascoli's theorem that for any
$\alpha'\in (0,\alpha)$, $v_R\in C^{0,\alpha'}(B_{R/2}(0)\cap\{x_1\leq
0\})$ and that, up to a subsequence,
\bequa\label{lim:co:ve}
\lim_{\eps\to 0}\eta_R\ve= v_R\hbox{ in
}C^{0,\alpha'}(B_{R/4}(0)\cap\{x_1\leq 0\})
\eequa
With (\ref{def:ve:sec3}) and (\ref{def:lae}), we have that
$(\eta_R\ve)(-\lae,\theta_\eps)=1$ for all $\eps>0$ and $R>0$ large
enough. Passing to the limit $\eps\to 0$ in this last equality, using
(\ref{lim:co:ve}) and (\ref{lim:lae}), we get that
$$v_R(-\lambda_0,\theta_0,0)=1$$
for $R>0$ large enough. With the same type of arguments, we get that $v\in
C^{0,\alpha}(\{x_1\leq 0\})$ and that $\lim_{R\to +\infty}v_R=v$ in
$C^{0,\alpha}_{loc}(\{x_1\leq 0\})$. Since $\eta_Rv=v_R$, we get that
$v(-\lambda_0,\theta_0)=1$. In particular, $v\not\equiv 0$ and
$\lambda_0>0$. This ends  Step 5.7.

\medskip\noindent{\bf  Step 5.8:} We claim that there exists $\theta\in
(0,1)$ such that $v\in C^{1,\theta}(\overline{\rnm})$ and
$$\ve\to v\hbox{ in }C^{1,\theta}_{loc}(\overline{\rnm})$$
when $\eps\to 0$. Indeed, it follows from  Step 5.7 that there exists
$\alpha>0$ such that for all $R>0$, there exists $C(R)>0$ such that
$$\Vert \ve\Vert_{C^{0,\alpha}(B_{R}(0)\cap\{x_1\leq 0\})}\leq C(R).$$
Following the proof of Proposition \ref{prop:app}, we let
$$\alpha_0:=\sup\{\alpha\in (0,1)/\, \forall R>0,\,\exists C(R)>0\hbox{
s.t. }\Vert \ve\Vert_{C^{0,\alpha}(B_{R}(0)\cap\{x_1\leq 0\})}\leq
C(R)\}.$$
We let $\alpha\in (0,\alpha_0)$ and $R>0$. We let $\tilde{R}>R$. There
exists $C(\tilde{R})>0$ such that
\bequa\label{ineq:ve:c1:1}
\Vert \ve\Vert_{C^{0,\alpha}(B_{\tilde{R}}(0)\cap\{x_1\leq 0\})}\leq
C(\tilde{R}).
\eequa
Since $\ve\equiv 0$ on $\partial\rnm$, we get with (\ref{ineq:ve:c1:1})
that
\bequa\label{ineq:est:123}
|\ve(x)|=|\ve(x)-\ve(x-(x_1,0))|\leq C(\tilde{R})|x_1|^\alpha
\eequa
for all $B_{\tilde{R}}(0)\cap\{x_1<0\}$ and all $\eps>0$. It then follows
from the properties of $\varphi$ (see (\ref{def:vphi}) with $x_0=0$) that
$$0\leq \fe(x):=\frac{(\eta\ve)^{\crit-1-\pe}}{\left|\frac{\varphi(\ke
x)}{\ke}\right|^s}\leq \frac{C}{\left|x\right|^{s-(\crit-1-\pe)\alpha}}$$
for all $\eps>0$ and all $x\in B_{\tilde{R}}(0)\cap\{x_1<0\}$. With the
properties (\ref{def:vphi}), we get that for any $\tilde{R}>0$ and any
$p>1$, we have that
$$\int_{B_{\tilde{R}}(0)\cap\{x_1<0\}}\frac{dx}{\left|\frac{\varphi(\ke
x)}{\ke}\right|^{p}}\leq C\int_{B_{\tilde{R}}(0)}\frac{dx}{|x|^p}$$
for all $\eps>0$ (note that the RHS can be infinite). Using the same
strategy as in the proof of Proposition \ref{prop:app}, we get that there
exists $\theta\in (0,1)$ such that  $v\in C^{1,\theta}(\overline{\rnm})$
and
$$\ve\to v\hbox{ in }C^{1,\theta}_{loc}(\overline{\rnm})$$
when $\eps\to 0$. We omit the proof and refer to the proof of Proposition
\ref{prop:app} for the details. This ends  Step 5.8.

\medskip\noindent{\bf  Step 5.9:} We claim that
\bequa\label{eq:lim:v}
\Delta v=\frac{v^{\crit-1}}{|x|^s}\hbox{ in }{\mathcal D}'(\rnm)
\eequa
and that
$$\int_{\rnm}|\nabla v|^2\,
dx=\mu_{s,a}(\Omega)^{\frac{\crit}{\crit-2}}=\mu_{s}(\rnm)^{\frac{\crit}{\crit-2}}.$$
Indeed, passing to the weak limit $\eps\to 0$ and then to the weak limit
$R\to +\infty$ in (\ref{eq:ve:sec3}), we get that
$$\Delta v=\frac{v^{\crit-1}}{|x|^s}\hbox{ in }{\mathcal D}'(\rnm).$$
Testing this equality with $v\in \hunrnm\setminus\{0\}$ and using the
optimal Hardy-Sobolev inequality (\ref{def:mus}), we get that
\bequa
\left(\int_{\rnm}|\nabla v|^2\,
dx\right)^{\frac{\crit-2}{\crit}}=\frac{\int_{\rnm}|\nabla v|^2\,
dx}{\left(\int_{\rnm}\frac{v^{\crit}}{|x|^s}\,
dx\right)^{\frac{2}{\crit}}}\geq\mu_{s}(\rnm).\label{ineq:sec3:rho}
\eequa
We then obtain that
\bequa\label{ineq:nrj:v:1}
\int_{\rnm}|\nabla v|^2\, dx\geq \mu_{s}(\rnm)^{\frac{\crit}{\crit-2}}.
\eequa
Since $0\leq\ve\leq 1$, it follows from Lebesgue's theorem that $\ve\to v$
strongly in $L^{\frac{2n}{n-2}}_{loc}(\rnm\cap\{x_1\leq 0\})$ when
$\eps\to 0$. Passing to the weak limit in (\ref{ineq:ve:1}) and using
(\ref{hyp:nrj:min}), we get that
\beqn
\int_{\rnm} |\nabla v_R|^2\, dx&\leq&
C_\alpha\Vert\nabla\eta_R\Vert_{n}^2\cdot\Vert v
\Vert_{L^{\frac{2n}{n-2}}(B_{2R}(0)\setminus B_{R}(0))}^{2}\nonumber\\
&&+(1+\alpha)\cdot(\lim_{\eps\to
0}\me^{\frac{\pe(n-2)}{\crit-2}})\mu_{s,a}(\Omega)^{\frac{\crit}{\crit-2}}\label{ineq:sec3:vr}
\eeqn
for all $R>0$. Since $v\in\hunrnm$, it follows from Sobolev's theorem that
$v\in L^{\frac{2n}{n-2}}(\rnm)$. Since
$\Vert\nabla\eta_R\Vert_{n}^2=\Vert\nabla\eta_1\Vert_{n}^2$ is independant
of $R>0$ and $v\in L^{\frac{2n}{n-2}}(\rnm)$, letting $R\to +\infty$ in
(\ref{ineq:sec3:vr}), we get that
\bequa\label{ineq:nrj:v:2}
\int_{\rnm} |\nabla v|^2\, dx\leq   (1+\alpha)\cdot(\lim_{\eps\to
0}\me^{\frac{\pe(n-2)}{\crit-2}})\mu_{s,a}(\Omega)^{\frac{\crit}{\crit-2}}
\eequa
Since $\alpha>0$ is arbitrary and $\me\leq 1$, we get with
(\ref{ineq:nrj:v:1}), (\ref{ineq:nrj:v:2}), Proposition \ref{prop:test}
and (\ref{ineq:sec3:rho}) that
$$\int_{\rnm}|\nabla v|^2\,
dx=\mu_{s,a}(\Omega)^{\frac{\crit}{\crit-2}}=\mu_{s}(\rnm)^{\frac{\crit}{\crit-2}},$$
and that
$$\lim_{\eps\to 0}\me^{\pe}=1.$$
This ends  Step 5.9. Proposition \ref{prop:sec3:3} then follows from Steps
5.1 to 5.9.
\end{proof}

\medskip\noindent{\bf   Step 5.10:} We claim that under the hypothesis of
Proposition \ref{prop:sec3:3}, we have that

\bequa\label{eq:limReps}
\lim_{R\to +\infty}\lim_{\eps\to 0}\int_{\Omega\setminus
B_{R\ke}(0)}\frac{\ue^{\crit-\pe}}{|x|^s}\, dx=0.
\eequa
Indeed, we let $R>0$. Since $D_0\varphi=Id_{\rn}$ and $\varphi(0)=0$, we
have that
$$\varphi\left(B_{\frac{R}{2}\ke}(0)\right)\subset B_{R\ke}(0)$$
for all $R>0$ and $\eps>0$ small enough. With a change of variable and
(\ref{hyp:nrj:min}), we get that
\beq
&&\int_{\Omega\setminus B_{R\ke}(0)}\frac{\ue^{\crit-\pe}}{|x|^s}\, dx\leq
\int_{\Omega\setminus
\varphi\left(B_{\frac{R}{2}\ke}(0)\right)}\frac{\ue^{\crit-\pe}}{|x|^s}\,
dx\\
&&\leq  \int_{\Omega}\frac{\ue^{\crit-\pe}}{|x|^s}\,
dx-\int_{\varphi\left(B_{\frac{R}{2}\ke}(0)\right)}\frac{\ue^{\crit-\pe}}{|x|^s}\,
dx\\
&&\leq
\mu_{s,a}(\Omega)^{\frac{\crit}{\crit-2}}+o(1)-\me^{-\pe\frac{(n-2)^2}{2(2-s)}}(1+o(1))\int_{B_{\frac{R}{2}}(0)}\frac{\ve^{\crit-\pe}}{|x|^s}\,
dx.
\eeq
when $\eps\to 0$. Letting $\eps\to 0$ and then $R\to +\infty$, we get with
(\ref{lim:re}) and Proposition \ref{prop:sec3:3} that
\beq
\lim_{R\to +\infty}\lim_{\eps\to 0}\int_{\Omega\setminus
B_{R\ke}(0)}\frac{\ue^{\crit-\pe}}{|x|^s}\, dx&\leq &
\mu_{s,a}(\Omega)^{\frac{\crit}{\crit-2}}-\lim_{R\to
+\infty}\int_{B_{\frac{R}{2}}(0)}\frac{v^{\crit}}{|x|^s}\, dx\\
&\leq &
\mu_{s,a}(\Omega)^{\frac{\crit}{\crit-2}}-\int_{\rnm}\frac{v^{\crit}}{|x|^s}\,
dx=0
\eeq
This last inequality yields (\ref{eq:limReps}).

\section{Refined Blow-Up analysis and strong pointwise estimates}
The objective of this section is the proof of the following strong
pointwise estimate
\begin{prop}\label{prop:fund:est}
Let $\Omega$ be a smooth bounded domain of $\rn$, $n\geq 3$. We let $s\in
(0,2)$. We let $(\pe)_{\eps>0}$ such that $\pe\in [0,\crit-2)$ for all
$\eps>0$ and (\ref{lim:pe}) holds. We consider $(\ue)_{\eps>0}\in \huno$
such that (\ref{hyp:ae}), (\ref{coerc:a}), (\ref{syst:ue}),
(\ref{hyp:nrj:min}) and (\ref{hyp:blowup}) hold. We let $\me$ as in
(\ref{def:me:xe}). Then, there exists $C>0$ such that
\bequa\label{eq:fund:est}
\ue(x)\leq C\cdot\left(\frac{\me}{\me^2+|x|^2}\right)^{\frac{n-2}{2}}
\eequa
for all $\eps>0$ and all $x\in\Omega$.
\end{prop}
This type of strong pointwise estimate first appeared in \cite{ha} in the
Euclidean context, and in \cite{hv1} in the Riemannian context. General
estimates are in \cite{dhr}.
\begin{proof} The rest of the section is mainly devoted to the proof of
the proposition. Here again, we follows the strategy of \cite{dhr}. We let
$(\ue)_{\eps>0}$ satisfying the hypothesis of Proposition
\ref{prop:fund:est}.

\medskip\noindent{\bf   Step 6.1:} We claim that there exists $C>0$ such
that
\bequa\label{ineq:est:1}
|x|^{\frac{n-2}{2}}\ue(x)^{1-\frac{\pe}{\crit-2}}\leq C
\eequa
for all $\eps>0$ and all $x\in\Omega$.

\noindent We proceed by contradiction and let $\ye\in\Omega$ such that
\bequa\label{hyp:wpe}
|\ye|^{\frac{n-2}{2}}\ue(\ye)^{1-\frac{\pe}{\crit-2}}=\sup_{x\in\Omega}|x|^{\frac{n-2}{2}}\ue(x)^{1-\frac{\pe}{\crit-2}}\to
+\infty
\eequa
when $\eps\to 0$. We let
\bequa\label{def:nueps}
\nu_\eps:=\ue(\ye)^{-\frac{2}{n-2}}\hbox{ and
}\ell_\eps:=\nu_\eps^{1-\frac{\pe}{\crit-2}}
\eequa
for all $\eps>0$. It follows from (\ref{hyp:wpe}) and (\ref{def:nueps})
that
\bequa\label{lim:ye:nueps}
\lim_{\eps\to 0}\nu_\eps=0\hbox{ and }\lim_{\eps\to
0}\frac{|\ye|}{\ell_\eps}=+\infty.
\eequa
It follows from (\ref{def:me:xe}) and (\ref{lim:re}) that
\bequa\label{lim:nu:pe}
\lim_{\eps\to 0}\nu_\eps^{\pe}=1.
\eequa
We let
\bequa\label{def:gamma:eps}
\gamma_{\eps}^{2}:=|y_{\eps}|^{s}|u_{\eps}(y_{\eps})|^{-(2^{*}-2-p_{\eps})},
\eequa
for all $\eps>0$. It follows from (\ref{lim:ye:nueps}) that
\bequa\label{lim:bis}
\lim_{\eps\to 0}\frac{\gamma_\eps}{|\ye|}=0.
\eequa

\smallskip\noindent{\it Case 6.1.1:} We assume that, up to a subsequence,
there exists $\rho>0$ such that
\bequa\label{hyp:case1:wpe}
\frac{d(\ye,\partial\Omega)}{\gamma_\eps}\geq 3\rho
\eequa
for all $\eps>0$. For any $x\in B_{2\rho}(0)$ and any $\eps>0$, we let
\bequa\label{def:we:wpe}
\we(x):=\nu_\eps^{\frac{n-2}{2}}\ue(\ye+\gamma_\eps x).
\eequa
Note that $\we$ is well defined thanks to (\ref{hyp:case1:wpe}). With
(\ref{hyp:wpe}) and (\ref{def:gamma:eps}), we get that
$$\left|\frac{\ye}{|\ye|}+\frac{\gamma_\eps}{|\ye|}x\right|^{\frac{n-2}{2}}\we(x)^{1-\frac{\pe}{\crit-2}}\leq
1.$$
In particular, with (\ref{lim:ye:nueps}), there exists $C_0>0$ such that
\bequa\label{bnd:we:wpe}
0\leq \we(x)\leq C_0
\eequa
for all $x\in B_{2\rho}(0)$ and all $\eps>0$. With (\ref{syst:ue}), we get
that
$$\Delta\we+\gamma_\eps^2 a_\eps(\ye+\gamma_\eps
x)\we=\frac{\we^{\crit-1-\pe}}{\left|\frac{\ye}{|\ye|}+\frac{\gamma_\eps}{|\ye|}x\right|^s}$$
for all $x\in B_{2\rho}(0)$ and all $\eps>0$. Since (\ref{lim:ye:nueps})
and (\ref{bnd:we:wpe}) hold, it follows from standard elliptic theory that
there exists $w\in C^1(B_{2\rho}(0))$ such that $w\geq 0$ and
\bequa\label{lim:we:wpe}
\lim_{\eps\to 0}\we=w
\eequa
in $C_{loc}^1(B_{2\rho}(0))$. It follows from (\ref{def:we:wpe}) that
$w(0)=1$. With a change of variable, we get that
\bequa\label{eq:int:wpe}
\int_{B_{\rho\gamma_\eps}(\ye)}\frac{\ue(x)^{\crit-\pe}}{|x|^s}\,
dx=\frac{\gamma_\eps^n\ue(\ye)^{\crit-\pe}}{|\ye|^s}\int_{B_{\rho}(0)}\frac{\we(x)^{\crit-\pe}}{\left|\frac{\ye}{|\ye|}+\frac{\gamma_\eps}{|\ye|}\cdot
x\right|^s}\, dx.
\eequa
With (\ref{def:gamma:eps}), (\ref{lim:nu:pe}), (\ref{lim:ye:nueps}) and
(\ref{def:nueps}), we then get that
$$\frac{\gamma_\eps^n\ue(\ye)^{\crit-\pe}}{|\ye|^s}=(1+o(1))\cdot
\left(\frac{|\ye|}{\ell_\eps}\right)^{\frac{s(n-2)}{2}}\to +\infty$$
when $\eps\to 0$. With (\ref{eq:int:wpe}), (\ref{lim:we:wpe}) and
(\ref{hyp:nrj:min}), we get that
$$\int_{B_{\rho}(0)}w^{\crit}\, dx=0,$$
and then $w\equiv 0$. A contradiction since $w(0)=1$. This ends Case
6.1.1.

\smallskip\noindent{\it Case 6.1.2:} We assume that
\bequa\label{hyp:case2:wpe}
\lim_{\eps\to 0}\frac{d(\ye,\partial\Omega)}{\gamma_\eps}=0.
\eequa
It then follows that there exists $y_0\in\partial\Omega$ such that
$$\lim_{\eps\to 0}\ye=y_0.$$
Since $y_0\in\partial\Omega$, which is smooth, we let $\varphi:U\to V$ as
in (\ref{def:vphi}) with $x_0=y_0$ and where $U,V$ are open neighborhoods
of $0$ and $y_0$ respectively. We let $\tue=\ue\circ\varphi$, which is
defined on $U\cap \{x_1\leq 0\}$. For any $i,j=1,...,n$, we let
$g_{ij}=(\partial_i\varphi,\partial_j\varphi)$, where $(\cdot,\cdot)$
denotes the Euclidean scalar product on $\rn$, and we consider $g$ as a
metric on $\rn$. We let $\Delta_g=-div_g(\nabla)$ the Laplace-Beltrami
operator with respect to the metric $g$. In the basis we choose,
$$\Delta_g=-g^{ij}\left(\partial_{ij}-\Gamma_{ij}^k\partial_k\right),$$
where $g^{ij}=(g^{-1})_{ij}$ are the coordinates of the inverse of the
tensor $g$ and the $\Gamma_{ij}^k$'s are the Christoffel symbols of the
metric $g$. As easily checked, we have that
$$\Delta_g\tue+\ae\circ\varphi(x)\cdot
\tue=\frac{\tue^{\crit-1-\pe}}{|\varphi(x)|^s}$$
weakly in $U\cap\{x_1<0\}$. We let $\ze\in\partial\Omega$ such that
\bequa\label{def:ze:bis}
|\ze-\ye|=d(\ye,\partial\Omega).
\eequa
We let $\tye,\tze\in U$ such that
\bequa\label{def:txe:tze:bis}
\varphi(\tye)=\ye\hbox{ and }\varphi(\tze)=\ze.
\eequa
It follows from the properties of $\varphi$ that
\bequa\label{ppty:txe:tze:bis}
\lim_{\eps\to 0}\tye=\lim_{\eps\to 0}\tze=0,\; (\tye)_1<0\hbox{ and
}(\tze)_1=0.
\eequa
At last, we let
$$\twe(x):=\frac{\tue(\tze+\gamma_\eps x)}{\tue(\tye)}$$
for all $x\in \frac{U-\tze}{\gamma_\eps}\cap \{x_1<0\}$. With
(\ref{ppty:txe:tze:bis}), we get that $\twe$ is defined on
$B_R(0)\cap\{x_1<0\}$ for all $R>0$, as soon as $\eps$ is small enough.
The function $\twe$ verifies
$$\Delta_{\tge}\twe+\gamma_\eps^2\ae\circ\varphi(\tze+\gamma_\eps
x)\twe=\frac{\twe^{\crit-1-\pe}}{\left|\frac{\varphi(\tze+\gamma_\eps
x)}{|y_{\eps}|}\right|^s}$$
weakly in $B_R(0)\cap\{x_1<0\}$. In this expression,
$\tge=g(\tze+\gamma_\eps x)$ and $\Delta_{\tge}$ is the Laplace-Beltrami
operator with respect to the metric $\tge$. With (\ref{hyp:case2:wpe}),
(\ref{def:ze:bis}) and (\ref{def:txe:tze:bis}), we get that
$$\varphi(\tze+\gamma_\eps x)=\ye+O_R(1)\gamma_\eps,$$
for all $x\in B_{R}(0)\cap\{x_1\leq 0\}$ and all $\eps>0$, where there
exists $C_R>0$ such that $|O_R(1)|\leq C_R$ for all $x\in
B_{R}(0)\cap\{x_1<0\}$. With (\ref{lim:bis}), we then get that
$$\lim_{\eps\to 0}\frac{|\varphi(\tze+\gamma_\eps x)|}{|\ye|}=1$$
in $C^0(B_{R}(0)\cap\{x_1\leq 0\})$. It then follows that
$$\Delta_{\tge}\twe+\gamma_{\eps}^2\ae\circ \varphi(\tze+\gamma_\eps
x)\twe=(1+o(1))\twe^{\crit-1-\pe}$$
weakly in $B_R(0)\cap\{x_1<0\}$, where $\lim_{\eps\to 0}o(1)=0$ in
$C^0(B_R(0)\cap\{x_1\leq 0\})$. Since $\twe$ vanishes on
$B_R(0)\cap\{x_1=0\}$ (in the sense of the trace) and that $0\leq \twe\leq
2$ (see for instance the proof of (\ref{bnd:we:wpe})), it follows from
standard elliptic theory that there exists $\tw\in C^1(B_R(0)\cap\{x_1\leq
0\})$ such that $$\lim_{\eps\to 0}\twe=\tw$$
in $C^{0}(B_{\frac{R}{2}}(0)\cap \{x_1\leq 0\})$. In particular,
\bequa\label{eq:v:vanish:bis}
\tw\equiv 0\hbox{ on }B_{\frac{R}{2}}(0)\cap \{x_1=0\}.
\eequa
Moreover, it follows from (\ref{hyp:case2:wpe}), (\ref{def:ze:bis}) and
(\ref{def:txe:tze:bis}) that
$$\twe\left(\frac{\tye-\tze}{\gamma_\eps}\right)=1\hbox{ and
}\lim_{\eps\to 0}\frac{\tye-\tze}{\gamma_\eps}=0.$$
In particular, $\tw(0)=1$. A contradiction with (\ref{eq:v:vanish:bis}).
This ends Case 6.1.2.

\medskip\noindent In both cases, we have contradicted (\ref{hyp:wpe}).
This proves (\ref{ineq:est:1}) and ends   Step 6.1.

\medskip\noindent As a remark, it follows from (\ref{syst:ue}),
(\ref{hyp:blowup}), (\ref{ineq:est:1}) and standard elliptic theory that
\bequa\label{cv:out:0}
\lim_{\eps\to 0}\ue=0\hbox{ in
}C^2_{loc}(\overline{\Omega}\setminus\{0\}).
\eequa
\medskip\noindent{\bf   Step 6.2:} This step is a slight improvement of
(\ref{ineq:est:1}). We claim that
\bequa\label{ineq:est:2}
\lim_{R\to +\infty}\lim_{\eps\to 0}\sup_{x\in \Omega\setminus
B_{R\ke}(0)}|x|^{\frac{n-2}{2}}\ue(x)^{1-\frac{\pe}{\crit-2}}=0.
\eequa
We proceed by contradiction and assume that there exists $\eps_0>0$ and a
family $(\ye)_{\eps>0}\in\Omega$ such that
\bequa\label{hyp:wpe:2}
|\ye|^{\frac{n-2}{2}}\ue(\ye)^{1-\frac{\pe}{\crit-2}}\geq\eps_0\hbox{ and
}\lim_{\eps\to 0}\frac{|\ye|}{\ke}=+\infty.
\eequa
We let
\bequa\label{def:nueps:2}
\nu_\eps:=\ue(\ye)^{-\frac{2}{n-2}}\hbox{ and
}\gamma_\eps:=\nu_\eps^{1-\frac{\pe}{\crit-2}}
\eequa
for all $\eps>0$. It follows from (\ref{cv:out:0}), (\ref{ineq:est:1}),
(\ref{hyp:wpe:2}) and (\ref{def:nueps:2}) that there exists
$\rho_0\in\mathbb{R}$ such that
\bequa\label{lim:ye:nueps:2}
\lim_{\eps\to 0}\ye=0,\;\lim_{\eps\to 0}\nu_\eps=0\hbox{ and
}\lim_{\eps\to 0}\frac{|\ye|}{\gamma_\eps}=\rho_0>0.
\eequa
Note that it follows from (\ref{def:me:xe}) and (\ref{lim:re}) that
\bequa\label{lim:nu:pe:2}
\lim_{\eps\to 0}\nu_\eps^{\pe}=1.
\eequa
We let $\varphi :U\to V$ as in (\ref{def:vphi}) with $x_0=0$ and where
$U,V$ are open neighborhoods of $0$. For any $x\in
\frac{U}{\gamma_\eps}\cap\{x_1<0\}$, we let
\bequa\label{def:we:wpe:2}
\bwe(x):=\nu_\eps^{\frac{n-2}{2}}\ue\circ\varphi(\gamma_\eps x).
\eequa
It follows from (\ref{ineq:est:1}) and the properties (\ref{def:vphi}) of
$\varphi$ that there exists $C>0$ such that
\bequa\label{bound:bwe}
|x|^{\frac{n-2}{2}}\bwe(x)^{1-\frac{\pe}{\crit-2}}\leq C
\eequa
for all $x\in \frac{U}{\gamma_\eps}\cap\{x_1<0\}$ and all $\eps>0$. As
above, we let the metric
$(\bge)_{ij}=(\partial_i\varphi,\partial_j\varphi)(\gamma_\eps x)$ for
$i,j=1,...,n$. With (\ref{syst:ue}), we get that
\bequa\label{eq:bwe}
\Delta_{\bge}\bwe+\gamma_\eps^2 a_\eps\circ\varphi(\gamma_\eps
x)\bwe=\frac{\bwe^{\crit-1-\pe}}{\left|\frac{\varphi(\gamma_\eps
x)}{\gamma_\eps}\right|^s}
\eequa
in $\frac{U}{\gamma_\eps}\cap\{x_1<0\}$ for all $\eps>0$. Moreover, $\bwe$
vanishes on $\frac{U}{\gamma_\eps}\cap\{x_1=0\}$. It then follows from
(\ref{bound:bwe}), (\ref{eq:bwe}) and standard elliptic theory (see for
instance \cite{gt}) that there exists $\bw\in
C^0(\rnm\cap\{x_1=0\})\setminus\{0\})$ such that $\bw\geq 0$ and
$$\lim_{\eps\to 0}\bwe=\bw$$
in $C^0(\rnm\cap\{x_1=0\})\setminus\{0\})$. We now write
$\ye=\varphi(\gamma_\eps\tye)$. It follows from (\ref{lim:ye:nueps:2})
that $\lim_{\eps\to 0}=y_0\neq 0$. As a consequence,
$$\bw(y_0)=\lim_{\eps\to 0}\bwe(\tye)=1,$$
and then $\bw\not\equiv 0$. We let $0<\delta<R$. With a change of
variable, we have that
\bequa\label{eq:int:wpe:2}
\lim_{\eps\to 0}\int_{\varphi(B_{R\gamma_\eps}(0))\setminus
\varphi(B_{\delta\gamma_\eps}(0))}\frac{\ue(x)^{\crit-\pe}}{|x|^s}\,
dx=\int_{B_{R}(0)\setminus B_{\delta}(0)}\frac{\bw(x)^{\crit}}{|x|^s}\,
dx.
\eequa
With (\ref{hyp:wpe:2}), we get that for any $\rho>0$, we have that
$$B_{\rho\ke}(0)\cap \left(\varphi(B_{R\gamma_\eps}(0))\setminus
\varphi(B_{\delta\gamma_\eps}(0))\right)=\emptyset$$
for all $\eps>0$ small enough, up to a subsequence. It then follows from
(\ref{eq:limReps}) that
$$\lim_{\eps\to 0}\int_{\varphi(B_{R\gamma_\eps}(0))\setminus
\varphi(B_{\delta\gamma_\eps}(0))}\frac{\ue(x)^{\crit-\pe}}{|x|^s}\,
dx=0.$$
This equality and (\ref{eq:int:wpe:2}) yield
$$\int_{B_{R}(0)\setminus B_{\delta}(0)}\frac{\bw^{\crit}}{|x|^s}\, dx=0$$
for all $R>\delta>0$. We then get that $\bw\equiv 0$. A contradiction
since $\bw(y_0)=1$. This ends   Step 6.2.

\medskip\noindent{\bf  Step 6.3:} We prove a first approximation of
(\ref{eq:fund:est}). More precisely, we claim that for any $\alpha\in
(0,n-2)$, there exists $C_\alpha>0$ such that
\bequa\label{ineq:est:nu}
|x|^{\alpha}\me^{\frac{n-2}{2}-\alpha}\ue(x)\leq C_\alpha
\eequa
for all $\eps>0$ and all $x\in\Omega$.
Indeed, since $\Delta+a$ is coercive on $\Omega$ and $(\ae)_{\eps>0}$
satisfies (\ref{hyp:ae}) and (\ref{coerc:a}), there exists $U_0$ an open
subset of $\rn$ such that $\overline{\Omega}\subset\subset U_0$, there
exists $\alpha_0>0$ and there exists $\lambda>0$ such that
\bequa\label{coerc:Ge:1}
\int_{U_0}\left(|\nabla\varphi|^2+(\ae-2\alpha_0)\varphi^2\right)\, dx\geq
\lambda \int_{U_0}\varphi^2\, dx
\eequa
for all $\varphi\in C^1_c(U_0)$ and all $\eps>0$. In other words, the
family of the operators $\Delta+\ae-\alpha_0$ is uniformly coercive in a
neighborhood of $\overline{\Omega}$. We let $G_\eps\in C^2(U_0\times
U_0\setminus\{(x,x)/x\in U_0\})$ be the Green's function for
$\Delta+\ae-\alpha_0$ with Dirichlet condition in $U_0$. In other words,
$G_\eps$ satisfies
\bequa\label{eq:Ge}
\Delta G_\eps(x,\cdot)+(\ae-\alpha_0)G_\eps(x,\cdot)=\delta_x
\eequa
weakly in ${\mathcal D}(U)$. It is standard that $G_\eps$ exists and,
since $0\in U$, that there exists $C>0$ such that
\bequa\label{ineq:Ge:1}
0<G_\eps(0,x)\leq C\cdot |x|^{2-n}
\eequa
for all $\eps>0$ and all $x\in \overline{U}\setminus\{0\}$. More
precisely, there exists $\delta_0>0$ and $C_0>0$ such that
\bequa\label{ineq:Ge:2}
G_\eps(0,x)\geq C_0\cdot |x|^{2-n}\hbox{ and }\frac{|\nabla
G_\eps(0,x)|}{|x|^{n-2}}\geq \frac{C_0}{|x|}
\eequa
for all $\eps>0$ and all $x\in B_{\delta_0}(0)\setminus\{0\}$. We let the
operator
$$L_\eps=\Delta +\left(\ae-\frac{\ue^{\crit-2-\pe}}{|x|^s}\right).$$

\medskip\noindent We claim that there exist $\nu_0\in (0,1)$ and $R_1>0$
such that for any $\nu\in (0,\nu_0)$ and any $R>R_1$, we have that
\bequa\label{ineq:LG}
L_\eps G_\eps^{1-\nu}>0
\eequa
for all $x\in \Omega\setminus B_{R\ke}(0)$ and for all $\eps>0$
sufficiently small. Indeed, we let $\nu_0\in (0,1)$ such that for any
$\nu\in (0,\nu_0)$, we have that
\bequa\label{def:nu}
\nu\cdot (\ae(x)-\alpha_0)\geq -\frac{1}{2}\alpha_0
\eequa
for all $\eps>0$ and all $x\in \Omega$. With (\ref{eq:Ge}), we get that
\bequa\label{eq:Ge:2}
\frac{L_\eps G_\eps^{1-\nu}}{G_\eps^{1-\nu}}(x)=\alpha_0+\nu\cdot
(\ae(x)-\alpha_0)+\nu\cdot(1-\nu)\cdot\frac{|\nabla
G_\eps|^2}{G_\eps^2}(x)-\frac{\ue(x)^{\crit-2-\pe}}{|x|^s}
\eequa
for all $x\in \Omega\setminus\{0\}$ and all $\eps>0$. It follows from the
pointwise estimate (\ref{ineq:est:2}) that there exists $R_1>0$ such that
for any $R>R_1$, we have that
\bequa\label{ineq:wpe:nu}
|x|^{2-s}\ue(x)^{\crit-2-\pe}\leq \frac{1}{2}\nu(1-\nu)C_0^2
\eequa
for all $\eps>0$ and all $x\in \Omega\setminus B_{R\ke}(0)$. Here, $C_0>0$
is as in (\ref{ineq:Ge:2}). We are now in position to prove
(\ref{ineq:LG}). We let $\nu\in(0,\nu_0)$ and $R>R_1$. We first let
$x\in\Omega$ such that $|x|\geq \delta_0$. It follows from (\ref{eq:Ge:2})
and (\ref{def:nu}) that
$$\frac{L_\eps G_\eps^{1-\nu}}{G_\eps^{1-\nu}}(x)\geq
\frac{\alpha_0}{2}-\frac{\ue(x)^{\crit-2-\pe}}{\delta_0^s}
$$
for all $\eps>0$. Inequality (\ref{ineq:LG}) then follows with this
inequality and (\ref{cv:out:0}). This proves (\ref{ineq:LG}) when
$|x|\geq\delta_0$.

\smallskip\noindent We let $x\in B_{\delta_0}(0)\setminus B_{R\ke}(0)$. It
follows from (\ref{eq:Ge:2}), (\ref{ineq:Ge:2}) and (\ref{ineq:wpe:nu})
that
$$\frac{L_\eps G_\eps^{1-\nu}}{G_\eps^{1-\nu}}(x)\geq
\frac{\alpha_0}{2}+\frac{\nu\cdot(1-\nu)\cdot
C_0^2}{|x|^2}-\frac{\nu\cdot(1-\nu)\cdot C_0^2}{2\cdot|x|^2}>0.$$
This proves (\ref{ineq:LG}) when $x\in B_{\delta_0}(0)\setminus
B_{R\ke}(0)$. Clearly these two assertions prove inequality
(\ref{ineq:LG}).

\medskip\noindent We let $R<R_1$ and $\nu\in (0,\nu_0)$. We claim that
there exists $C(R)>0$ such that
\bequa\label{ineq:LG:2}\left\{\begin{array}{ll}
L_\eps\left(C(R)\me^{\frac{n-2}{2}-\nu(n-2)}G_\eps(0,\cdot)^{1-\nu}\right)>L_\eps\ue&\hbox{
in }\Omega\setminus B_{R\ke}(0)\\
&\\
C(R)\me^{\frac{n-2}{2}-\nu(n-2)}G_\eps(0,\cdot)^{1-\nu}>\ue&\hbox{ on
}\partial\Omega\setminus B_{R\ke}(0))
\end{array}\right\}
\eequa
Indeed, the first inequality is trivial since $L_\eps\ue=0$ and
(\ref{ineq:LG}) holds. Concerning the second inequality, we get with the
definition (\ref{def:me:xe}) of $\me$, the limit (\ref{lim:re}) and
(\ref{ineq:Ge:2}) that
\beq
\frac{\ue(x)}{\me^{\frac{n-2}{2}-\nu(n-2)}G_\eps(0,x)^{1-\nu}}&\leq
&C_0^{\nu-1}\cdot\me^{-(n-2)(1-\nu)}\cdot|x|^{(n-2)(1-\nu)}\\
&&\leq 2\cdot C_0^{1-\nu}\cdot R^{(n-2)(1-\nu)}:=C(R)
\eeq
for all $x\in\Omega\cap \partial B_{R\ke}(0)$. The inequalities
(\ref{ineq:LG:2}) are proved.

\medskip\noindent Since $G_\eps(0,x)^{1-\nu}>0$ in $\overline{\Omega\cap
\partial B_{R\ke}(0)}$ and $L_\eps G_\eps(0,x)^{1-\nu}>0$ in $\Omega\cap
\partial B_{R\ke}(0)$, it follows from \cite{bnv} that $L_\eps$ verifies
the comparison principle. It then follows from (\ref{ineq:LG:2}) that
$$\ue(x)\leq C(R)\me^{\frac{n-2}{2}-\nu(n-2)}G_\eps(0,x)^{1-\nu}$$
for all $x\in \Omega\setminus \overline{B}_{R\ke}(0)$. With
(\ref{ineq:Ge:1}), we get that there exists $C'(R)>0$ such that
$$\ue(x)\leq C'(R)\me^{\frac{n-2}{2}-\nu(n-2)}|x|^{2-n+\nu(n-2)}$$
for all $x\in \Omega\setminus \overline{B}_{R\ke}(0)$. Up to taking
$C'(R)$ larger, it follows from (\ref{def:me:xe}) that this inequality
holds on the whole set $\Omega$. Taking $\alpha=(n-2)\cdot(1-\nu)$, we get
(\ref{ineq:est:nu}) for $\alpha$ close to $n-2$. As easily checked, this
implies the inequality for all $\alpha\in (0,n-2)$. This ends the proof of
(\ref{ineq:est:nu}).

\medskip\noindent{\bf  Step 6.4:} We are in position to prove Proposition
\ref{prop:fund:est}. For all $\eps>0$, we let $\ye\in \Omega$ such that
$$\max_{x\in\Omega}|x|^{n-2}\ue(\xe)\ue(x)=|\ye|^{n-2}\ue(\xe)\ue(\ye).$$
Clearly, Proposition \ref{prop:fund:est} is equivalent to proving that
\bequa\label{lim:maxue}
|\ye|^{n-2}\ue(\xe)\ue(\ye)=O(1)
\eequa
when $\eps\to 0$.

\medskip\noindent{\it Case 6.4.1:} We assume that
$$|\ye|=O(\ke)$$
when $\eps\to 0$. We then get with (\ref{def:me:xe}) that
$$|\ye|^{n-2}\ue(\xe)\ue(\ye)=O(1).$$
when $\eps\to 0$. This proves (\ref{lim:maxue}) in Case 6.4.1.

\medskip\noindent{\it Case 6.4.2:} We assume that
\bequa\label{lim:case542}
\lim_{\eps\to 0}\frac{|\ye|}{\ke}=+\infty.
\eequa
As in the beginning of  Step 6.3, we choose $U_0$ such that
$\overline{\Omega}\subset\subset U_0$ such that $\Delta+\ae$ is coercive
on $U_0$. We let $H_\eps$ be the Green's function for $\Delta+\ae$ on
$U_0$ with Dirichlet boundary condition. It follows from Green's
representation formula and standard estimates on the Green's function that
\bequa\label{green:1}
\ue(x)\leq \int_\Omega
H_\eps(x,y)\cdot\frac{\ue(y)^{\crit-1-\pe}}{|y|^s}\, dy\leq C\int_\Omega
|x-y|^{2-n}\cdot\frac{\ue(y)^{\crit-1-\pe}}{|y|^s}\, dy
\eequa
for all $x\in\Omega$. We let
$$\hat{v}_\eps(x)=\me^{\frac{n-2}{2}}\ue(\ke x)$$
for all $x\in \ke^{-1}\Omega$ and all $\eps>0$. It follows from
Proposition \ref{prop:sec3:3} and (\ref{ineq:est:nu}) that for any
$\alpha\in (0,n-2)$, there exists $C_\alpha>0$ such that
\bequa\label{ef:ve}
\hat{v}_\eps(x)\leq \frac{C_\alpha}{1+|x|^\alpha}
\eequa
for all $x\in \ke^{-1}\Omega$ and all $\eps>0$. It follows from
(\ref{green:1}) and a change of variable that
\beqn
\me^{-\frac{n-2}{2}}\ue(\ye)&\leq & C \int_{\ke^{-1}\Omega} |\ye-\ke
y|^{2-n}\frac{\hat{v}_\eps(y)^{\crit-1-\pe}}{|y|^s}\, dx\nonumber\\
&\leq & C \int_{\ke^{-1}\Omega\cap \left\{|\ye-\ke
y|\geq\frac{|\ye|}{2}\right\}} \frac{1}{|\ye-\ke y|^{n-2}}\cdot
\frac{\hat{v}_\eps(y)^{\crit-1-\pe}}{|y|^s}\, dx\nonumber\\
&& +C \int_{\ke^{-1}\Omega\cap \left\{|\ye-\ke y|<\frac{|\ye|}{2}\right\}}
\frac{1}{|\ye-\ke y|^{n-2}}\cdot
\frac{\hat{v}_\eps(y)^{\crit-1-\pe}}{|y|^s}\, dx.\label{ef:split}
\eeqn
We estimate the two integrals of the RHS separately. With (\ref{ef:ve}),
we get that
\beqn
&&\int_{\ke^{-1}\Omega\cap \{|\ye-\ke y|\geq\frac{|\ye|}{2}\}}
\frac{1}{|\ye-\ke y|^{n-2}}\cdot
\frac{\hat{v}_\eps(y)^{\crit-1-\pe}}{|y|^s}\, dx\nonumber\\
&&\leq C\cdot |\ye|^{2-n}\int_{\ke^{-1}\Omega}
\frac{1}{|y|^s(1+|y|^{\alpha\cdot(\crit-1-\pe)})}\, dy\nonumber\\
&&\leq C\cdot |\ye|^{2-n}\label{ef:split:1}
\eeqn
for all $\eps>0$ small and $\alpha$ close enough to $n-2$. On the other
hand, with (\ref{ef:ve}), we get that
\beq
&&\int_{\ke^{-1}\Omega\cap \{|\ye-\ke y|\leq\frac{|\ye|}{2}\}}
\frac{1}{|\ye-\ke y|^{n-2}}\cdot
\frac{\hat{v}_\eps(y)^{\crit-1-\pe}}{|y|^s}\, dx\\
&&\leq C\int_{\ke^{-1}\Omega\cap \{|\ye-\ke y|\leq\frac{|\ye|}{2}\}}
\frac{1}{|\ye-\ke y|^{n-2}}\cdot \frac{1}{|y|^{\alpha(\crit-1-\pe)+s}}\,
dx\\
&&\leq \frac{C\cdot
\ke^{\alpha(\crit-1-\pe)+s}}{|\ye|^{\alpha(\crit-1-\pe)+s}}\int_{\ke^{-1}\Omega\cap
\{|\ye-\ke y|\leq\frac{|\ye|}{2}\}} \frac{1}{|\ye-\ke y|^{n-2}}\, dy\\
&&\leq \frac{C\cdot
\ke^{\alpha(\crit-1-\pe)+s}}{|\ye|^{\alpha(\crit-1-\pe)+s}}\cdot
\frac{|\ye|^2}{|\ke|^n}\\
&&\leq C |\ye|^{2-n}\cdot
\left(\frac{\ke}{|\ye|}\right)^{(\crit-1-\pe)\alpha+s-n}.
\eeq
Since $\lim_{\alpha\to n-2}\lim_{\eps\to 0}(\crit-1-\pe)\alpha+s-n=2-s>0$,
we get with (\ref{lim:case542}) and $\alpha$ close enough to $n-2$ that
\bequa\label{ef:split:2}
\int_{\ke^{-1}\Omega\cap \{|\ye-\ke y|\geq\frac{|\ye|}{2}\}}
\frac{1}{|\ye-\ke y|^{n-2}}\cdot
\frac{\hat{v}_\eps(y)^{\crit-1-\pe}}{|y|^s}\,
dx=o\left(|\ye|^{2-n}\right),
\eequa
when $\eps\to 0$. Plugging together (\ref{ef:split:1}) and
(\ref{ef:split:2}) into (\ref{ef:split}), we get that
$$\me^{-\frac{n-2}{2}}\ue(\ye)=O\left(|\ye|^{2-n}\right)$$
when $\eps\to 0$. This proves that (\ref{lim:maxue}) holds in Case 6.4.2.

\medskip\noindent In both cases, we have proved that (\ref{lim:maxue})
holds. As easily checked, (\ref{eq:fund:est}) and then Proposition
\ref{prop:fund:est} follow from (\ref{lim:maxue}) and (\ref{def:me:xe}).
This ends  Step 6.4, and therefore proves Proposition
\ref{prop:fund:est}.\end{proof}

\medskip\noindent{\bf  Step 6.5:} From Proposition \ref{prop:fund:est}, we
can derive pointwise estimates for $\ve$. This is the object of the
following proposition
\begin{prop}\label{prop:est:c1}
Assume that the hypothesis of Proposition \ref{prop:fund:est} are
satisfied. Then there exists $C>0$ such that
$$\ve(x)\leq \frac{C}{(1+|x|^2)^{\frac{n-2}{2}}}\hbox{ and
}|\nabla\ve|(x)\leq \frac{C}{(1+|x|^2)^{\frac{n-1}{2}}}$$
for all $\eps>0$ and all $x\in \frac{U}{\ke}\cap \{x_1<0\}$, where $\ve$
was defined in (\ref{def:ve:sec3}) and $U$ is as in (\ref{def:vphi}) with
$x_0=0$.
\end{prop}
\begin{proof} The first inequality of the proposition is an immediate
consequence of the estimate (\ref{eq:fund:est}) and the definition
(\ref{def:ve:sec3}) of $\ve$. Concerning the second inequality, we proceed
by contradiction and assume that there exists a family $(\ye)_{\eps>0}$
such that $\ye\in U$ for all $\eps\to 0$ and such that
\bequa\label{hyp:lim:c1:infty}
\lim_{\eps\to
0}\left(1+\left|\frac{\ye}{\ke}\right|\right)^{n-1}\left|\nabla\ve\left(\frac{\ye}{\ke}\right)\right|=+\infty.
\eequa

\medskip\noindent{\it Case 6.5.1:} we assume that $\ye\not\to 0$ when
$\eps\to 0$. It follows from the pointwise estimate (\ref{eq:fund:est})
that for any $\delta>0$, there exists $C(\delta)>0$ such that
$$\ue(x)\leq C(\delta)\me^{\frac{n-2}{2}}$$
for all $x\in \overline{\Omega}\setminus B_{\delta}(x_0)$ and all
$\eps>0$. We then get that
$$\Delta (\me^{\frac{2-n}{2}}\ue)+\ae \cdot (\me^{\frac{2-n}{2}}\ue)
=\me^{\frac{n-2}{2}(\crit-2-\pe)}\frac{(\me^{\frac{2-n}{2}}\ue)^{\crit-1-\pe}}{|x|^s}$$
in ${\mathcal D}'(\Omega\setminus \bar{B}_{\delta}(x_0))$. It then follows
from standard elliptic theory that
\bequa\label{eq:c1:1}
\Vert\me^{\frac{2-n}{2}}\ue\Vert_{C^{2}(\overline{\Omega}\setminus
B_{3\delta}(x_0))}=O(1)
\eequa
when $\eps\to 0$. Since $\ye\not\to 0$, there exists $\delta>0$ such that,
up to a subsequence, $|\ye|\geq 4\delta$ for $\eps>0$. It follows from
(\ref{eq:c1:1}) that $\nabla\ue(\varphi(\ye))=O(\me^{\frac{n-2}{2}})$ when
$\eps\to 0$. A contradiction with (\ref{hyp:lim:c1:infty}). This proves
the Proposition in Case 6.5.1.

\medskip\noindent{\it Case 6.5.2:} We assume that
\bequa\label{hyp:Case2}
\lim_{\eps\to 0}\ye=0\hbox{ and }\lim_{\eps\to
0}\frac{|\ye|}{\ke}=+\infty.
\eequa
We let $\varphi$ as in (\ref{def:vphi}) with $x_0=0$ and define
$$\he(x):=\frac{|\ye|^{n-2}}{\ke^{\frac{n-2}{2}}}\ue\circ\varphi(|\ye|
x)$$
for all $x\in \frac{U}{|\ye|}\cap \{x_1\leq 0\}$. It follows from
(\ref{eq:fund:est}) and (\ref{lim:re}) that there exists $C>0$ such that
\bequa\label{upper:he}
\he(x)\leq C\cdot |x|^{2-n}
\eequa
for all $x\in \frac{U}{|\ye|}\cap \{x_1\leq 0\}$, $x\neq 0$. We let
$$\Delta_{\bar{g}_\eps}=\bar{g}_\eps^{ij}\left(\partial_{ij}-\Gamma_{ij}^k(\bar{g}_\eps)\partial_k\right),$$
the Laplace-Beltrami operator for the metric
$(\bar{g}_\eps)_{ij}=(\partial_i\varphi,\partial_j\varphi)(\ke x)$. In
this expression, the $\bar{g}_\eps^{ij}=(\bar{g}_\eps^{-1})_{ij}$ are the
coordinates of the inverse of the tensor $\bar{g}_\eps $ and the
$\Gamma_{ij}^k(\bar{g}_\eps)$ are the Christoffel symbols associated to
the metric $\bar{g}_\eps$. After a change of variables, (\ref{syst:ue})
rewrites as
$$\Delta_{\bar{g}_\eps}\he+|\ye|^2\ae(\varphi(|\ye|x))\he=\ke^{\pe\frac{n-2}{2}}\left(\frac{\ke}{|\ye|}\right)^{2-s-\pe(n-2)}\frac{\he^{\crit-1-\pe}}{\left|\frac{\varphi(|\ye|
x)}{|\ye|}\right|^s}$$
in ${\mathcal D}'\left(\frac{U}{|\ye|}\cap \{x_1< 0\}\right)$. Since
(\ref{lim:re}), (\ref{hyp:Case2}) and (\ref{upper:he}) hold and since
$s\in (0,2)$, there exists $p>\frac{n}{2}$ such that
$$\Delta_{\bar{g}_\eps}\he+|\ye|^2\ae(\varphi(|\ye|x))\he=\fe\hbox{ in
}{\mathcal D}'\left(\frac{U}{|\ye|}\cap \{x_1< 0\}\right),$$
where $\fe\in L^p_{loc}(\frac{U}{|\ye|}\cap \{x_1\leq 0\}\setminus\{0\})$
uniformly wrt $\eps\to 0$. Since $\he\equiv 0$ on $\frac{U}{|\ye|}\cap
\{x_1= 0\}$ and (\ref{upper:he}) holds, it follows from standard elliptic
theory that there exists for any $\delta_1>\delta_2>0$, there exists
$C'(\delta_1,\delta_2)>0$ such that
$$\Vert\he\Vert_{C^{1}((B_{\delta_1}(0)\setminus
B_{\delta_2}(0))\cap\{x_1\leq 0\})}\leq C'(\delta_1,\delta_2)$$
for all $\eps>0$. It then follows that
$$\left\vert\nabla\he\left(\frac{\ye}{|\ye|}\right)\right\vert=O(1)$$
when $\eps\to 0$. Coming back to the definitions of $\he$ and $\ve$, we
get a contradiction with (\ref{hyp:lim:c1:infty}). This proves the
Proposition in Case 6.5.2.

\medskip\noindent{\it Case 6.5.3:} We assume that
$$|\ye|=O(\ke)$$
when $\eps\to 0$. In this case, It follows from Proposition
\ref{prop:sec3:3} that
$\left|\nabla\ve\left(\frac{\ye}{\ke}\right)\right|=O(1)$ when $\eps\to
0$. We get a contradiction with (\ref{hyp:lim:c1:infty}). This proves the
Proposition in Case 3.

\medskip\noindent In all the cases, we have contradicted
(\ref{hyp:lim:c1:infty}). This proves Proposition \ref{prop:est:c1}.
\end{proof}
\begin{coro}\label{coro:lim:0} Let $(\ue)_{\eps>0}$ as in the hypothesis
of Proposition \ref{prop:fund:est}. Then there exists $H\in 
C^1(\overline{\Omega}\setminus\{0\})$ such that
$$\ue(\xe)\ue\to H\hbox{ in }C^1_{loc}(\overline{\Omega}\setminus\{0\})$$
when $\eps\to 0$.
\end{coro}
\begin{proof} We let $H_\eps(x):=\ue(\xe)\ue(x)$ for all $x\in\Omega$ and 
all $\eps>0$. It follows from Proposition \ref{prop:fund:est} that for any 
open subset $U$ such that $\overline{U}\subset 
\overline{\Omega}\setminus\{0\}$, there exists $C(U)>0$ such that 
$|H_\eps(x)|\leq C(U)$ for all $x\in U$ and all $\eps>0$. Equation 
(\ref{syst:ue}) rewrites as
$$\Delta H_\eps+\ae 
H_\eps=\ue(\xe)^{2+\pe-\crit}\frac{H_\eps^{\crit-1-\pe}}{|x|^s}$$
in $\Omega$. The conclusion of the Corollary is then a consequence of 
standard elliptic theory.
\end{proof}

\section{Pohozaev identity and proof of the theorems}
In this section, we prove the following estimate:
\begin{prop}\label{prop:poho}
Let $\Omega$ be a smooth bounded domain of $\rn$, $n\geq 4$. We let $s\in 
(0,2)$. We let $(\pe)_{\eps>0}$ such that $\pe\in [0,\crit-2)$ for all 
$\eps>0$ and (\ref{lim:pe}) holds. We consider $(\ue)_{\eps>0}\in \huno$ 
such that (\ref{hyp:ae}), (\ref{coerc:a}), (\ref{syst:ue}), 
(\ref{hyp:nrj:min}) and (\ref{hyp:blowup}) hold. We let $\me$ as in 
(\ref{def:me:xe}) and $v$ as in Proposition \ref{prop:sec3:3}. Then, we 
have that
$$\lim_{\eps\to 
0}\frac{\pe}{\me}=\frac{(n-s)\int_{\partial\rnm}|x|^2|\nabla v|^2\, 
dx}{n(n-2)^2\mu_{s}(\rnm)^{\frac{n-s}{2-s}}}\cdot H(0).$$
In this expression, $H(0)$ is the mean curvature of the oriented boundary 
$\partial\Omega$ at $0$.
\end{prop}
We prove the Proposition in the sequel, and postpone the proofs of 
Theorems \ref{th:intro} and \ref{th:eq} to the end of the section. We let 
$\pe\geq 0$ such that $\lim_{\eps\to 0}\pe=0$. We let $\ue$, $\ae$
and $a$ as in (\ref{hyp:ae}), (\ref{coerc:a}), (\ref{syst:ue}),
(\ref{hyp:nrj:min}) and (\ref{hyp:blowup}). We assume that $0<s<2$ and let
$\xe$, $\me$, $\ke$ as in (\ref{def:me:xe}). Since $\lim_{\eps\to
0}\xe=0$, we consider the chart $\varphi$ defined in (\ref{def:vphi}) with
$x_0=0$.

\medskip\noindent{\bf  Step 7.1:} We provide a Pohozaev-type identity for
$\ue$. It follows from Proposition \ref{prop:app} that $\ue\in
C^1(\overline{\Omega})$ and that $\Delta\ue\in L^p(\Omega)$ for all $p\in
(1,\frac{n}{s})$. In the sequel, we denote by $\nu(x)$ the outward normal
vector at $x\in\partial\Omega$ of the oriented hypersurface
$\partial\Omega$ (oriented as the boundary of $\Omega$). Integrating by
parts, we get that
\beq
&&\int_\Omega x^i\partial_i\ue\Delta\ue\, dx\\
&&=-\int_{\partial \Omega}x^i\partial_i\ue\partial_{\nu}\ue\,
d\sigma+\int_{\Omega}\partial_j(x^i\partial_i\ue)\partial_j\ue\, dx\\
&&=-\int_{\partial \Omega}x^i\partial_i\ue\partial_{\nu}\ue\,
d\sigma+\int_{\Omega}|\nabla\ue|^2\, dx+\int_\Omega
x^i\partial_i\frac{|\nabla\ue|^2}{2}\, dx\\
&&=\left(1-\frac{n}{2}\right)\int_{\Omega}|\nabla\ue|^2\,
dx+\int_{\partial
\Omega}\left((x,\nu)\frac{|\nabla\ue|^2}{2}-x^i\partial_i\ue\partial_{\nu}\ue\right)\,
d\sigma\\
&&=\left(1-\frac{n}{2}\right)\left(\int_{\partial
\Omega}\ue\partial_{\nu}\ue\, d\sigma+\int_{\Omega}\ue\Delta\ue\,
dx\right)\\
&&+\int_{\partial
\Omega}\left((x,\nu)\frac{|\nabla\ue|^2}{2}-x^i\partial_i\ue\partial_{\nu}\ue\right)\,
d\sigma.
\eeq
Using the equation (\ref{syst:ue}) in the RHS, we get that
\beqn
&&\int_\Omega x^i\partial_i\ue\Delta\ue\,
dx=\left(1-\frac{n}{2}\right)\left(\int_{\Omega}\frac{\ue^{\crit-\pe}}{|x|^s}\,
dx-\int_{\Omega}\ae\ue^2\, dx\right)\nonumber\\
&&+\int_{\partial
\Omega}\left(\left(1-\frac{n}{2}\right)\ue\partial_{\nu}\ue+(x,\nu)\frac{|\nabla\ue|^2}{2}-x^i\partial_i\ue\partial_{\nu}\ue\right)\,
d\sigma.\label{eq:poho:1}
\eeqn
On the other hand, using the equation (\ref{syst:ue}) satisfied by $\ue$,
we get that
\beqn
&&\int_\Omega x^i\partial_i\ue\Delta\ue\, dx=\int_\Omega
x^i\partial_i\ue\frac{\ue^{\crit-1-\eps}}{|x|^s}\, dx-\int_\Omega
x^i\partial_i\ue\ae\ue\, dx\nonumber\\
&&=\int_\Omega
x^i|x|^{-s}\partial_i\left(\frac{\ue^{\crit-\pe}}{\crit-\pe}\right)\,
dx-\int_\Omega x^i\partial_i\ue\ae\ue\, dx\nonumber\\
&& =-\int_\Omega\partial_i(x^i|x|^{-s})\frac{\ue^{\crit-\pe}}{\crit-\pe}\,
dx+\int_{\partial \Omega}
\frac{(x,\nu)}{\crit-\pe}\cdot\frac{\ue^{\crit-\pe}}{|x|^s}\,
d\sigma-\int_\Omega x^i\partial_i\ue\ae\ue\, dx\nonumber\\
&& =-\int_\Omega\frac{n-s}{|x|^s}\cdot\frac{\ue^{\crit-\pe}}{\crit-\pe}\,
dx+\frac{1}{2}\int_\Omega (n\ae +x^i\partial_i\ae)\ue^2\, dx\nonumber\\
&&+\int_{\partial \Omega}
\frac{(x,\nu)}{\crit-\pe}\cdot\frac{\ue^{\crit-\pe}}{|x|^s}\,
d\sigma-\int_{\partial \Omega}\frac{(x,\nu)}{2}\ae\ue^2 \,
d\sigma.\label{eq:poho:2}
\eeqn
Plugging together (\ref{eq:poho:1}) and (\ref{eq:poho:2}), we get that
\beqn
&&\left(\frac{n-2}{2}-\frac{n-s}{\crit-\pe}\right)\int_{\Omega}\frac{\ue^{\crit-\pe}}{|x|^s}\,
dx+\int_\Omega\left(\ae+\frac{(x,\nabla\ae)}{2}\right)\ue^2\,
dx\nonumber\\
&&=\int_{\partial
\Omega}\left(-\frac{n-2}{2}\ue\partial_\nu\ue+(x,\nu)\frac{|\nabla\ue|^2}{2}\right.\nonumber\\
&&\left.-x^i\partial_i\ue\partial_{\nu}\ue-\frac{(x,\nu)}{\crit-\pe}\cdot\frac{\ue^{\crit-\pe}}{|x|^s}\right)\,
d\sigma+\int_{\partial \Omega}\frac{(x,\nu)}{2}\ae\ue^2\,
dx\label{eq:poho:3}
\eeqn
for all $\eps>0$. Since $\ue\equiv 0$ on $\partial\Omega$, we get that
\beqn
&&\frac{(n-2)\pe}{2\cdot(\crit-\pe)}\int_{\Omega}\frac{\ue^{\crit-\pe}}{|x|^s}\,
dx-\int_\Omega\left(\ae+\frac{(x,\nabla\ae)}{2}\right)\ue^2\,
dx\nonumber\\
&&=\frac{1}{2}\int_{\partial \Omega}(x,\nu)|\nabla\ue|^2\,
d\sigma\label{eq:poho:5}.
\eeqn

\medskip\noindent{\bf  Step 7.2:} We first deal with the RHS of
(\ref{eq:poho:5}). We take $\varphi$ as in (\ref{def:vphi}) with $x_0=0$.
With the pointwise limit of Corollary \ref{coro:lim:0}, we get that
$$\int_{\partial \Omega}(x,\nu)|\nabla\ue|^2\, d\sigma=\int_{\partial
\Omega\cap\varphi(U)}(x,\nu)|\nabla\ue|^2\, d\sigma +o(\me)$$
when $\eps\to 0$ as soon as $n\geq 4$. With a change of variable, we get 
that
\beqn
&&\int_{\partial \Omega}(x,\nu)|\nabla\ue|^2\, d\sigma=\nonumber\\
&&(1+o(1))\cdot \int_{D_\eps}\left(\frac{\varphi(\ke
x)}{\ke},\nu\circ\varphi(\ke x)\right)|\nabla\ve|_{\tge}^2\sqrt{|\tge|}\,
dx\nonumber\\
&&+o(\me^{n-2})\label{eq:poho:rhs:1}
\eeqn
where the metric $\tge$ is such that
$(\tge)_{ij}=(\partial_i\varphi,\partial_j\varphi)(\ke x)$ for all
$i,j=2,...,n$, $|\tge|=\det(\tge)$ and
$$D_\eps=\frac{U}{\ke}\cap\{x_1=0\}.$$
Using the expression of $\varphi$ (see (\ref{def:vphi})), we get that
$$\nu(\varphi(x))=\frac{(1,-\partial_2\varphi_0(x), ...,
-\partial_n\varphi_0(x))}{\sqrt{1+\sum_{i=2}^n(\partial_i\varphi_0(x))^2}}$$
for all $x\in U\cap\{x_1=0\}$. We then get that
$$(\nu\circ\varphi(x),\vec{X})=(1+O(|x|^2))\cdot\left(X^1-\sum_{i=2}^nX^i\partial_i\varphi_0(x)\right)$$
for all $x\in U\cap\{x_1=0\}$ and all $\vec{X}\in\rn$. In this expression
$O(1)$ is bounded for $x\in U\cap\{x_1=0\}$ and $\vec{X}\in\rn$. With the
expression of $\varphi$ (see (\ref{def:vphi})), we get that
\beqn
&&(\varphi(\ke x),\nu\circ\varphi(\ke x))\nonumber\\
&&=(1+O(\ke^2|x|^2))\left(\varphi_0(\ke
x)-\ke\sum_{i=2}^nx^i\partial_i\varphi_0(\ke x)\right)\nonumber\\
&&=(1+O(\ke^2|x|^2))\cdot\left(-\frac{1}{2}\ke^2\partial_{ij}\varphi(0)x^ix^j+O(1)(\ke^3|x|^3)\right)\label{dev:PS}
\eeqn
for $\eps>0$ and $x\in \frac{U}{\ke}\cap\{x_1=0\}$. Plugging
(\ref{dev:PS}) into (\ref{eq:poho:rhs:1}), using the estimates of
Proposition \ref{prop:fund:est}, Lebesgue's convergence theorem and
letting $\eps\to 0$, we get that
\bequa\label{eq:poho:rhs:2}
\int_{\partial \Omega}(x,\nu)|\nabla\ue|^2\,
d\sigma=\left(-\frac{1}{2}\int_{\partial\rnm}\partial_{ij}\varphi_0(0)x^ix^j|\nabla
v|^2\, dx+o(1)\right)\cdot\ke
\eequa
when $n\geq 4$ and where $\lim_{\eps\to 0}o(1)=0$.

\medskip\noindent{\bf  Step 7.3:} It follows from Proposition
\ref{prop:fund:est}
\bequa\label{norm:L2}
\int_\Omega\ue^2\, dx=o(\me)
\eequa
when $\eps\to 0$ and as soon as $n\geq 4$. Plugging (\ref{eq:poho:rhs:2}) 
into
(\ref{eq:poho:5}), using (\ref{hyp:nrj:min}) and (\ref{norm:L2}), we get
that
\beqn
&&\left(\frac{n-2}{2\cdot\crit}\mu_{s}(\rnm)^{\frac{n-s}{2-s}}+o(1)\right)\pe\nonumber\\
&&=\left(-\frac{1}{4}\int_{\partial\rnm}\partial_{ij}\varphi_0(0)x^ix^j|\nabla
v|^2\, dx+o(1)\right)\cdot\me\label{eq:poho:7}
\eeqn
where $\lim_{\eps\to 0}o(1)=0$ and when $n\geq 4$. With (\ref{eq:poho:7}),
we get that
\beq
&&\lim_{\eps\to
0}\frac{n-2}{2\cdot\crit}\mu_{s}(\rnm)^{\frac{n-s}{2-s}}\cdot\frac{\pe}{\me}\\
&&=-\frac{1}{4}\int_{\partial\rnm}\partial_{ij}\varphi_0(0)x^ix^j|\nabla
v|^2\, dx
\eeq
when $n\geq 4$. We consider the second fondamental form associated to
$\partial\Omega$, namely
$$II_p(x,y)=(d\nu_px,y)$$
for all $p\in\partial\Omega$ and all $x,y\in T_{p}\partial\Omega$ (recall
that $\nu$ is the outward normal vector at the hypersurface
$\partial\Omega$). In the canonical basis of
$\partial\rnm=T_0\partial\Omega$, the matrix of the bilinear form $II_{0}$
is $-D^2_0\varphi_0$, where $D^2_0\varphi_0$ is the Hessian matrix of
$\varphi_0$ at $0$. With this remark and (\ref{eq:poho:7}), we get that
\bequa\label{lim:pe:ke:0}
\lim_{\eps\to
0}\frac{\pe}{\me}=\frac{(n-s)}{(n-2)^2}\mu_{s}(\rnm)^{-\frac{n-s}{2-s}}\cdot\int_{\partial\rnm}II_{0}(x,x)|\nabla
v|^2\, dx
\eequa
when $n\geq 4$. Since $v\geq 0$, that $v\in C^2(\rnm)$ and $v$ verifies 
(\ref{eq:lim:v}), it follows from the strong maximum principle that $v>0$ 
in $\rnm$. Moreover, it follows from the definition (\ref{def:ve:sec3}) 
and the pointwise estimate (\ref{eq:fund:est}) that there exists $C>0$ 
such that
$$v(x)\leq\frac{C}{(1+|x|^2)^{\frac{n-2}{2}}}$$
for all $x\in\rnm$. We let 
$\tilde{v}(x):=|x|^{2-n}v\left(\frac{x}{|x|^2}\right)$ be the Kelvin 
transform of $v$. As easily checked, $\tilde{v}\in 
C^2(\overline{\rnm}\setminus\{0\})$ and verifies
$$\Delta \tilde{v}=\frac{\tilde{v}^{\crit-1}}{|x|^s}\hbox{ in }\rnm\hbox{ 
and }\tilde{v}(x)\leq\frac{C}{(1+|x|^2)^{\frac{n-2}{2}}}$$
for all $x\in\rnm$. Since $\tilde{v}$ vanishes on $\partial\rnm$, it then 
follows from standard elliptic theory that $\tilde{v}\in 
C^1(\overline{\rnm})$ and then, that there exists $C>0$ such that 
$\tilde{v}(x)\leq C|x|$ for all $x\in B_1(0)\cap\rnm$. Coming back to the 
function $v$, we get that tehre exists $C>0$ such that
$$v(x)\leq\frac{C}{(1+|x|^2)^{\frac{n-1}{2}}}$$
for all $x\in\rnm$. It follows from Proposition \ref{prop:sym} of Appendix 
B that there exists $w\in C^2(\rr_{-}^\star\times \rr)$ such that 
$v(x_1,x')=w(x_1,|x'|)$ for all $(x_1,x')\in \rr_-^\star\times \rr^{n-1}$. 
In particular, $|\nabla v|(0,x')$ is radially symmetrical wrt 
$x'\in\partial\rnm$. Since we have chosen a chart $\varphi$ that is 
Euclidean at $0$, we get that
\beq
\int_{\partial\rnm}II_{0}(x,x)|\nabla v|^2\,
dx&=&\frac{\sum_{i=2}^n(II_0)^{ii}}{n}\int_{\partial\rnm}|x|^2|\nabla 
v|^2\,
dx\\
&=&\frac{H(0)}{n}\int_{\partial\rnm}|x|^2|\nabla v|^2\,
dx.
\eeq
Note that we have used here that in the chart $\varphi$ defined in 
(\ref{def:vphi}), the matrix of the first fundamental form at $0$ is the 
identity. Plugging thsi last inequality in (\ref{lim:pe:ke:0}), we get 
that
\bequa\label{lim:pe:ke}
\lim_{\eps\to 0}\frac{\pe}{\me}=\frac{(n-s)\int_{\partial\rnm}|x|^2|\nabla 
v|^2\, dx}{n(n-2)^2\mu_{s}(\rnm)^{\frac{n-s}{2-s}}}\cdot H(0)
\eequa
when $n\geq 4$.

\medskip\noindent{\bf  Step 7.4:} We are now in position to prove Theorems 
\ref{th:intro} and \ref{th:eq}. We prove Theorem \ref{th:intro} by 
contradiction and assume that there are no extremals for (\ref{def:mus}). 
It follows from Propositions \ref{prop:subcrit} and \ref{prop:min:nonzero} 
that there exists $\ue\in\huno$ such that (\ref{syst:ue}), 
(\ref{hyp:nrj:min}) and (\ref{hyp:blowup}) hold with $\ae\equiv 0$ and 
$\pe=\eps$. Since $0<s<2$, then (\ref{lim:pe:ke}) holds with $\pe=\eps$ 
when $n\geq 4$. We then get that $H(0)\geq 0$. A contradiction with the 
assumptions of Theorem \ref{th:intro}. This proves the first point of 
Theorem \ref{th:intro} when $n\geq 4$. Concerning the compactness, any 
sequence of minimizers of (\ref{def:mus}) satisfies (\ref{syst:ue}) and 
(\ref{hyp:nrj:min}) with $\pe\equiv 0$ and $a\equiv 0$. If the sequence of 
minimizers blows up, we get with (\ref{lim:pe:ke}) that $H(0)=0$. A 
contradiction with our initial assumption. Then we get that the sequence 
does not blow up. It then follows from standard elliptic theory that it 
converges in $\huno$. This proves Theorem \ref{th:intro} when $n\geq 4$.

\medskip\noindent Concerning Theorem \ref{th:eq}, the proof is quite 
similar to that of Theorem \ref{th:intro}. We assume that the conclusion 
of the theorem does not hold. It follows from Propositions 
\ref{prop:subcrit} and \ref{prop:min:nonzero} that there exists 
$\ue\in\huno$ such that (\ref{syst:ue}), (\ref{hyp:nrj:min}) and 
(\ref{hyp:blowup}) hold with $\ae\equiv a$ and $\pe=\eps$. The proof is 
then the same as the proof of Theorem \ref{th:intro}.

\section{Appendix: Regularity of weak solutions}
In this appendix, we prove the following regularity result:
\begin{prop}\label{prop:app}
Let $\Omega$ be a smooth bounded domain of $\rn$, $n\geq 3$. We let $s\in
(0,2)$ and $a\in C^0(\overline{\Omega})$. We let $\eps\in [0,\crit-2)$ and
consider $u\in\huno$ a weak solution of
$$\Delta u+au=\frac{|u|^{\crit-2-\eps}u}{|x|^s}  \hbox{ in }{\mathcal
D}'(\Omega).$$
Then there exists $\theta\in (0,1)$ such that $u\in
C^{1,\theta}(\overline{\Omega})$.
\end{prop}
\begin{proof}

\medskip\noindent{\bf  Step 8.1:} We follow the strategy developed by
Trudinger. Let $\beta\geq 1$, and $L>0$. We let
$$G_L(t)=\left\{\begin{array}{ll}
|t|^{\beta-1}t& \hbox{ if }|t|\leq L\\
\beta L^{\beta-1}(t-L)+L^\beta & \hbox{ if }t\geq L\\
\beta L^{\beta-1}(t+L)-L^\beta & \hbox{ if }t\leq -L
\end{array}\right.$$
and
$$H_L(t)=\left\{\begin{array}{ll}
|t|^{\frac{\beta-1}{2}}t& \hbox{ if }|t|\leq L\\
\frac{\beta+1}{2}L^{\frac{\beta-1}{2}}(t-L)+L^{\frac{\beta+1}{2}} & \hbox{
if }t\geq L\\
\frac{\beta+1}{2}L^{\frac{\beta-1}{2}}(t+L)-L^{\frac{\beta+1}{2}} & \hbox{
if }t\leq -L
\end{array}\right.$$
As easily checked,
$$0\leq t G_L(t)\leq H_L(t)^2\hbox{ and
}G_L'(t)=\frac{4\beta}{(\beta+1)^2}(H_L'(t))^2$$
for all $t\in\rr$ and all $L>0$.
Let $\eta\in C^\infty_c(\rn)$. As easily checked, $\eta^2G_L(u),\eta
H_L(u)\in\huno$. With the equation verified by $u$, we get that
\bequa\label{eq:C}
\int_\Omega\nabla u\nabla (\eta^2 G_L(u))\, dx=\int_\Omega
\frac{|u|^{\crit-2-\eps}}{|x|^s}\eta^2 u G_L(u)\, dx-\int_\Omega a\eta^2
uG_L(u)\, dx.
\eequa
We let $J_L(t)=\int_0^t G_L(\tau)\, d\tau$ for all $t\in\rr$. Integrating
by parts, we get that
\beqn
&&\int_\Omega\nabla u\nabla (\eta^2 G_L(u))\, dx=\int_\Omega\eta^2
G_L'(u)|\nabla u|^2\, dx+\int_\Omega \nabla\eta^2\nabla J_L(u)\,
dx\nonumber\\
&&=\frac{4\beta}{(\beta+1)^2}\int_\Omega \eta^2 |\nabla H_L(u)|^2\,
dx+\int_\Omega(\Delta\eta^2)J_L(u)\, dx\nonumber\\
&&= \frac{4\beta}{(\beta+1)^2}\int_\Omega |\nabla (\eta H_L(u))|^2\,
dx+\frac{4\beta}{(\beta+1)^2}\int_\Omega\eta\Delta\eta |H_L(u)|^2\,
dx\nonumber\\
&&+\int_\Omega(\Delta\eta^2)J_L(u)\, dx\label{eq:A}
\eeqn
On the other hand, with H\"older's inequality and the definition of
$\mu_{s}(\rn)$, we then get that
\beqn
&&\int_\Omega \left(\frac{|u|^{\crit-2-\eps}}{|x|^s}-a\right)\cdot \eta^2
u G_L(u)\, dx\leq \int_\Omega\left( |a|+
\frac{|u|^{\crit-2-\eps}}{|x|^s}\right)\cdot(\eta H_L(u))^2\,
dx\nonumber\\
&&\leq \left(\int_{\Omega\cap\hbox{Supp }\eta}\frac{(|a|\cdot
|x|^s+|u|^{\crit-2-\pe})^{\frac{\crit-\eps}{\crit-2-\eps}}}{|x|^s}\right)^{1-\frac{2}{\crit-\eps}}\nonumber\\
&&\times \left(\int_{\Omega}\frac{|\eta
H_L(u)|^{\crit}}{|x|^s}\right)^{\frac{2}{\crit}}\times\left(\int_\Omega\frac{dx}{|x|^s}\right)^{\frac{2\eps}{\crit\cdot(\crit-\eps)}}\nonumber\\
&&\leq \alpha\cdot \int_{\Omega}|\nabla(\eta H_L(u))|^2\, dx\label{eq:B}
\eeqn
where
\begin{eqnarray*}
\alpha&:=&\left(\int_{\Omega\cap\hbox{Supp }\eta}\frac{(|a|\cdot
|x|^s+|u|^{\crit-2-\pe})^{\frac{\crit-\eps}{\crit-2-\eps}}}{|x|^s}\,
dx\right)^{1-\frac{2}{\crit-\eps}}\\
&&\times\mu_{s}(\rn)^{-1}\left(\int_\Omega\frac{dx}{|x|^s}\right)^{\frac{2\eps}{\crit\cdot(\crit-\eps)}}
\end{eqnarray*}
Plugging (\ref{eq:A}) and (\ref{eq:B}) into (\ref{eq:C}), we get that
\bequa\label{eq:D}
A\cdot \int_\Omega |\nabla (\eta H_L(u))|^2\, dx\leq
\frac{4\beta}{(\beta+1)^2}\int_\Omega|\eta\Delta\eta| |H_L(u)|^2\,
dx+\int_\Omega|\Delta(\eta^2)J_L(u)|\, dx
\eequa
where
\begin{eqnarray*}
A&:=&\frac{4\beta}{(\beta+1)^2} -\left(\int_{\Omega\cap\hbox{Supp
}\eta}\frac{(|a|\cdot
|x|^s+|u|^{\crit-2-\pe})^{\frac{\crit-\eps}{\crit-2-\eps}}}{|x|^s}\,
dx\right)^{1-\frac{2}{\crit-\eps}}\\
&&\times
\mu_{s}(\rn)^{-1}\left(\int_\Omega\frac{dx}{|x|^s}\right)^{\frac{2\eps}{\crit\cdot(\crit-\eps)}}
\end{eqnarray*}

\medskip\noindent{\bf  Step 8.2:} We let $p_0=\sup\{p\geq 1/\, u\in
L^p(\Omega)\}$. It follows from Sobolev's embedding theorem that
$p_0\geq\frac{2n}{n-2}$. We claim that
$$p_0=+\infty.$$
We proceed by contradiction and assume that
$$p_0<\infty.$$
Let $p\in (2,p_0)$. It follows from the definition of $p_0$ that $u\in
L^p(\Omega)$. Let $\beta=p-1>1$. For any $x\in\overline{\Omega}$, we let
$\delta_x>0$ such that
\begin{eqnarray}
&&\left(\int_{\Omega\cap B_{2\delta_x}(x)}\frac{(|a|\cdot
|x|^s+|u|^{\crit-2-\pe})^{\frac{\crit-\eps}{\crit-2-\eps}}}{|x|^s}\,
dx\right)^{1-\frac{2}{\crit-\eps}}\mu_{s}(\rn)^{-1}\nonumber\\
&&\times
\left(\int_\Omega\frac{dx}{|x|^s}\right)^{\frac{2\eps}{\crit\cdot(\crit-\eps)}}\leq
\frac{2\beta}{(\beta+1)^2}.\label{eq:E}
\end{eqnarray}
Since $\overline{\Omega}$ is compact, we get that there exists
$x_1,...,x_N\in\overline{\Omega}$ such that
$$\overline{\Omega}\subset \bigcup_{i=1}^N B_{\delta_{x_i}}(x_i).$$
We fix $i\in\{1,...,N\}$ and let $\eta\in
C^\infty(B_{2\delta_{x_i}}(x_i))$ such that $\eta(x)=1$ for all $x\in
B_{\delta_{x_i}}(x_i)$. We then get with (\ref{eq:D}) and (\ref{eq:E})
that
\beqn
&&\frac{2\beta}{(\beta+1)^2}\int_\Omega |\nabla (\eta H_L(u))|^2\,
dx\nonumber\\
&&\leq \frac{4\beta}{(\beta+1)^2}\int_\Omega|\eta\Delta\eta| |H_L(u)|^2\,
dx+\int_\Omega|\Delta\eta^2|\cdot |J_L(u)|\, dx.\label{eq:F}
\eeqn
Recall that it follows from Sobolev's inequality that there exists
$K(n,2)>0$ that depends only on $n$ such that
\bequa\label{ineq:sob}
\left(\int_{\rn}|f|^{\frac{2n}{n-2}}\, dx\right)^{\frac{n-2}{n}}\leq
K(n,2)\int_{\rn}|\nabla f|^2\, dx
\eequa
for all $f\in \hunrn$. It follows from (\ref{eq:F}) and (\ref{ineq:sob})
that
\beq
&&\frac{2\beta}{(\beta+1)^2}K(n,2)^{-1} \left(\int_\Omega |\eta
H_L(u)|^{\frac{2n}{n-2}}\, dx\right)^{\frac{n-2}{n}}\\
&&\leq \frac{4\beta}{(\beta+1)^2}\int_\Omega|\eta\Delta\eta| |H_L(u)|^2\,
dx+\int_\Omega|\Delta\eta^2|\cdot |J_L(u)|\, dx
\eeq
for all $L>0$. As easily checked, there exists $C_0>0$ such that
$|J_L(t)|\leq C_0\cdot |t|^{\beta+1}$ for all $t\in\rr$ and all $L>0$.
Since $u\in L^{\beta+1}(\Omega)$, we get that there exists a constant
$C=C(\eta,u,\beta,\Omega)$ independant of $L$ such that
$$\int_{\Omega\cap B_{\delta_{x_i}}(x_i)} |H_L(u)|^{\frac{2n}{n-2}}\,
dx\leq \int_\Omega |\eta H_L(u))|^{\frac{2n}{n-2}}\, dx\leq C$$
for all $L>0$. Letting $L\to +\infty$, we get that
$$\int_{\Omega\cap B_{\delta_{x_i}}(x_i)} |u|^{\frac{n}{n-2}(\beta+1)}\,
dx<+\infty,$$
for all $i=1...N$. We then get that $u\in
L^{\frac{n}{n-2}(\beta+1)}(\Omega)=L^{\frac{n}{n-2}p}(\Omega)$. And then,
$\frac{n}{n-2}p\leq p_0$ for all $p\in (2,p_0)$. Letting $p\to p_0$, we
get a contradiction. Then $p_0=+\infty$ and $u\in L^p(\Omega)$ for all
$p\geq 1$. This ends  Step 8.2.

\medskip\noindent{\bf  Step 8.3:} We claim that
$$u\in C^{0,\alpha}(\overline{\Omega})$$
for all $\alpha\in (0,1)$. Indeed, it follows from  Step 8.2 and the
assumption $0<s<2$ that there exists $p>\frac{n}{2}$ such that
$$\fe:=\frac{|u|^{\crit-2-\eps}u}{|x|^s}-au\in L^p(\Omega).$$
It follows from standard elliptic theory that, in this case, $u\in
C^{0,\alpha}(\overline{\Omega})$ for all $\alpha\in (0,\min\{2-s,1\})$. We
let
$$\alpha_0=\sup\{\alpha\in (0,1)/\, u\in
C^{0,\alpha}(\overline{\Omega})\}.$$
We let $\alpha\in (0,\alpha_0)$. Then $u\in
C^{0,\alpha}(\overline{\Omega})$. Since $u(0)=0$, we then get that
\bequa\label{ineq:app:1}
|u(x)|\leq |u(x)-u(0)|\leq C |x|^\alpha.
\eequa
We then get with (\ref{ineq:app:1}) that
$$\left|\fe(x)\right|=\left|\frac{|u(x)|^{\crit-1-\eps}u}{|x|^s}-au\right|\leq
\frac{C}{|x|^{s-(\crit-1-\eps)\alpha}}$$
for all $x\in\Omega$. We distinguish 2 cases:

\smallskip\noindent{\it Case 8.3.1:} $s-(\crit-1-\eps)\alpha_0\leq 0$. In
this case, for any $p>1$, up to taking $\alpha$ close enough to
$\alpha_0$, we get that
$$\fe\in L^p(\Omega).$$
Since $\Delta u+au=\fe$ and $u\in\huno$, it follows from standard elliptic
theory that there exist exists $\theta\in (0,1)$ such that $u\in
C^{1,\theta}(\overline{\Omega})$. It follows that $\alpha_0=1$. This
proves the claim in Case 8.3.1.

\smallskip\noindent{\it Case 8.3.2:} $s-(\crit-1-\eps)\alpha_0>0$. In this
case, for any $p<\frac{n}{s-(\crit-1-\eps)\alpha_0}$, up to taking
$\alpha$ close enough to $\alpha_0$, we get that
$$\fe\in L^p(\Omega).$$
We distinguish 3 subcases.

\smallskip\noindent{\it Case 8.3.2.1:} $s-(\crit-1-\eps)\alpha_0<1$. In
this case, up to taking $\alpha$ close enough to $\alpha_0$, there exists
$p>n$ such that
$$\fe\in L^p(\Omega).$$
Since $\Delta u=\fe$ and $u\in\huno$, it follows from standard elliptic
theory that there exist exists $\theta\in (0,1)$ such that $u\in
C^{1,\theta}(\overline{\Omega})$. It follows that $\alpha_0=1$. This
proves the claim in Case 8.3.2.1.

\smallskip\noindent{\it Case 8.3.2.2:} $s-(\crit-1-\eps)\alpha_0=1$. In
this case, for any $p<n$, up to taking $\alpha$ close enough to
$\alpha_0$, we get that
$$\fe\in L^p(\Omega).$$
Since $\Delta u+au=\fe$ and $u\in\huno$, it follows from standard elliptic
theory that $u\in C^{0,\tilde{\alpha}}(\overline{\Omega})$ for all
$\tilde{\alpha}\in (0,1)$. It follows that $\alpha_0=1$. This proves the
claim in Case 8.3.2.2.

\smallskip\noindent{\it Case 8.3.2.3:}  $s-(\crit-1-\eps)\alpha_0>1$. In
this case, it follows from standard elliptic theory that $u\in
C^{0,\tilde{\alpha}}(\overline{\Omega})$ for all
$$\tilde{\alpha}\leq 2-(s-(\crit-1-\eps)\alpha_0).$$
It follows from the definition of $\alpha_0$ that
$$\alpha_0\geq  2-(s-(\crit-1-\eps)\alpha_0),$$
and then
$$0\geq 2-s+\left(\crit-2-\eps\right)\alpha_0>0,$$
a contradiction since $s<2$ and $\eps<\crit-2$. This proves that Case
7.3.2.3 does not occur, and we are back to the other cases.

\smallskip\noindent Clearly, theses cases end  Step 8.3.

\medskip\noindent{\bf  Step 8.4:} We claim that there exists $\theta\in
(0,1)$ such that
$$u\in C^{1,\theta}(\overline{\Omega}).$$
We proceed as in  Step 8.3. We let $\alpha\in (0,1)$ (note that
$\alpha_0=1$). We then get that
$$\left|\fe(x)\right|=\left|\frac{|u(x)|^{\crit-1-\eps}u}{|x|^s}\right|\leq
\frac{C}{|x|^{s-(\crit-1-\eps)\alpha}}$$
for all $x\in\Omega$. We distinguish 2 cases:

\smallskip\noindent{\it Case 8.4.1:} $s-(\crit-1-\eps)\leq 0$. In this
case, for any $p>1$, up to taking $\alpha$ close enough to $\alpha_0$, we
get that
$$\fe\in L^p(\Omega).$$
Since $\Delta u+au=\fe$ and $u\in\huno$, it follows from standard elliptic
theory that there exist exists $\theta\in (0,1)$ such that $u\in
C^{1,\theta}(\overline{\Omega})$. It follows that $\alpha_0=1$. This
proves the claim in Case 8.4.1.

\smallskip\noindent{\it Case 8.4.2:} $s-(\crit-1-\eps)>0$. In this case,
for any $p<\frac{n}{s-(\crit-1-\eps)}$, up to taking $\alpha$ close enough
to $1$, we get that
$$\fe\in L^p(\Omega).$$
As easily checked,
$$1-(s-(\crit-1-\eps))=2-s+(\crit-1-\eps)-1>\crit-2-\eps.$$
We the get that there exists $p>n$ such that $\fe\in L^p(\Omega)$. Since
$\Delta u+au=\fe$ and $u\in\huno$, it follows from standard elliptic
theory that there exists $\theta\in (0,1)$ such that $u\in
C^{1,\theta}(\overline{\Omega})$. This proves the claim in Case 8.4.2.

\medskip\noindent Combining Case 8.4.1 and Case 8.4.2, we obtain  Step 
8.4.
Proposition \ref{prop:app} then follows from Step 8.4.
\end{proof}

\end{document}